\documentclass[preprint,12pt]{elsarticle}
\usepackage{amsmath, amssymb, amsthm}
\usepackage{color}
\usepackage{latexsym}
\usepackage{indentfirst}
\usepackage{graphicx}
\usepackage{placeins}
\usepackage{booktabs}
\usepackage{algorithm}
\usepackage{algorithmic}
\usepackage{multirow}

\theoremstyle{plain}
\newtheorem{theorem}{Theorem}

\theoremstyle{definition}
\newtheorem{defn}[theorem]{Definition}

\theoremstyle{remark}
\newtheorem*{remark}{Remark}

\newcommand{\etal}{\textit{et al}{}}

\newcommand{\RR}{\mathbb{R}}

\newcommand{\Or}{\mathcal{O}}
\newcommand{\mI}{\mathcal{I}}
\newcommand{\mJ}{\mathcal{J}}
\newcommand{\mP}{\mathcal{P}}
\newcommand{\mC}{\mathcal{C}}
\newcommand{\mD}{\mathcal{D}}
\newcommand{\TT}{\mathrm{T}}
\newcommand{\veps}{\varepsilon}

\newcommand{\abs}[1]{\lvert#1\rvert}
\newcommand{\norm}[1]{\lVert#1\rVert}

\newcommand{\bvec}[1]{\ensuremath{\mathbf{#1}}}

\newcommand{\MATLAB}{\textsf{MATLAB}}
\newcommand{\NL}{\mathsf{NL}}
\newcommand{\IL}{\mathsf{IL}}
\newcommand{\mc}[1]{\mathcal{#1}}
\newcommand{\wh}[1]{\widehat{#1}}
\newcommand{\wt}[1]{\widetilde{#1}}

\journal{Journal of Computational Physics}

\begin{document}

\begin{frontmatter}



\title{Fast construction of
  hierarchical matrix representation from matrix-vector
  multiplication}


\author[ll]{Lin Lin} 
\ead{linlin@math.princeton.edu}

\author[jl]{Jianfeng Lu} 
\ead{jianfeng@cims.nyu.edu}

\author[ly]{Lexing Ying}
\ead{lexing@math.utexas.edu}

\address[ll]{Program in Applied and Computational Mathematics, Princeton
  University, Princeton, NJ 08544. }

  \address[jl]{Department of Mathematics, Courant Institute of Mathematical
  Sciences, New York University, 251 Mercer St., New York, NY 10012.  }

\address[ly]{Department of Mathematics and ICES, University of Texas at
  Austin, 1 University Station, C1200, Austin, TX 78712. }

\begin{abstract}

  We develop a hierarchical matrix construction algorithm using
  matrix-vector multiplications, based on the randomized singular
  value decomposition of low-rank matrices. The algorithm uses
  $\Or(\log n)$ applications of the matrix on structured random test
  vectors and $\Or(n \log n)$ extra computational cost, where $n$ is
  the dimension of the unknown matrix. Numerical examples on
  constructing Green's functions for elliptic operators in two
  dimensions show efficiency and accuracy of the proposed algorithm.
  

\end{abstract}


\end{frontmatter}



\section{Introduction}


In this work, we consider the following problem: Assume that an
unknown symmetric matrix $G$ has the structure of a hierarchical
matrix ($\mc{H}$-matrix) \cite{BormGrasedyckHackbusch:06,
  Hackbusch:99, HackbuschKhoromskijSauter:00}, that is, certain
off-diagonal blocks of $G$ are low-rank or approximately low-rank (see
the definitions in Sections~\ref{sec:1Dalg} and
\ref{sec:Hmatrix}). The task is to construct $G$ efficiently only from
a ``black box'' matrix-vector multiplication subroutine (which shall
be referred to as \textsf{matvec} in the following). In a slightly
more general setting when $G$ is not symmetric, the task is to
construct $G$ from ``black box'' matrix-vector multiplication
subroutines of both $G$ and $G^{\TT}$. In this paper, we focus on the
case of a symmetric matrix $G$. The proposed algorithm can be extended to the
non-symmetric case in a straightforward way.

\subsection{Motivation and applications}

Our motivation is mainly the situation that $G$ is given as the
Green's function of an elliptic equation. In this case, it is proved
that $G$ is an $\mc{H}$-matrix under mild regularity assumptions
\cite{BebendorfHackbusch:03}. For elliptic equations, methods like
preconditioned conjugate gradient, geometric and algebraic multigrid
methods, sparse direct methods provide application of the matrix $G$
on vectors. The algorithm proposed in this work then provides an
efficient way to construct the matrix $G$ explicitly in the
$\mc{H}$-matrix form.

Once we obtain the matrix $G$ as an $\mc{H}$-matrix, it is possible to
apply $G$ on vectors efficiently, since the application of an
$\mc{H}$-matrix on a vector is linear scaling. Of course, for elliptic
equations, it might be more efficient to use available fast solvers
directly to solve the equation, especially if only a few right hand
sides are to be solved. However, sometimes, it would be advantageous
to obtain $G$ since it is then possible to further compress $G$
according to the structure of the data (the vectors that $G$ will be
acting on), for example as in numerical homogenization
\cite{EngquistRunborg:02}. Another scenario is that the data has
special structure like sparsity in the choice of basis, the
application of the resulting compressed matrix will be more efficient
than the ``black box'' elliptic solver.

Let us remark that, in the case of elliptic equations, it is also
possible to use the $\mc{H}$-matrix algebra to invert the direct
matrix (which is an $\mc{H}$-matrix in \textit{e.g.} finite element
discretization).  Our method, on the other hand, provides an efficient
alternative algorithm when a fast matrix-vector multiplication is
readily available.  From a computational point of view, what is
probably more attractive is that our algorithm facilitates a
parallelized construction of the $\mc{H}$-matrix, while the direct
inversion has a sequential nature \cite{Hackbusch:99}.

As another motivation, the purpose of the algorithm is to recover the
matrix via a ``black box'' matrix-vector multiplication subroutine. A
general question of this kind will be that under which assumptions of
the matrix, one can recover the matrix efficiently by matrix-vector
multiplications.  If the unknown matrix is low-rank, the recently
developed randomized singular value decomposition algorithms
\cite{HalkoMartinssonTropp:09,
  LibertyWoolfeMartinssonRokhlinTygert:07,
  WoolfeLibertyRokhlinTygert:08} provide an efficient way to obtain
the low-rank approximation through application of the matrix on random
vectors. Low-rank matrices play an important role in
many applications.  However, the assumption is too strong in
many cases that the whole matrix is low-rank. Since the class of $\mc{H}$-matrices
is a natural generalization of the one of low-rank matrices, the
proposed algorithm can be viewed as a further step in this direction.

\subsection{Randomized singular value decomposition algorithm}

A repeatedly leveraged tool in the proposed algorithm is the
randomized singular value decomposition algorithm for computing a low
rank approximation of a given numerically low-rank matrix. This has
been an active research topic in the past several years with vast
literature. For the purpose of this work, we have adopted the
algorithm developed in \cite{LibertyWoolfeMartinssonRokhlinTygert:07},
although other variants of this algorithm with similar ideas can also
be used here. For a given matrix $A$ that is numerically low-rank,
this algorithm goes as following to compute a rank-$r$ factorization.

\begin{algorithm}[ht]
\begin{small}
\begin{flushleft}
\begin{minipage}{4.8in}
\begin{algorithmic}[1]
  \STATE Choose a Gaussian random matrix $R_1 \in \RR^{n\times
    (r+c)}$ where $c$ is a small constant; 
  \STATE Form $A R_1$ and apply SVD to $A R_1$. The first $r$
  left singular vectors give $U_1$;
  \STATE Choose a Gaussian random matrix $R_2 \in \RR^{n\times
    (r +c)}$;
  \STATE Form $R_2^{\TT} A$ and apply SVD to $A^{\TT} R_2$. The first 
  $r$ left singular vectors give $U_2$;
  \STATE $M = (R_2^{\TT}U_1)^{\dagger} [R_2^{\TT}(A R_1)] (U_2^{\TT} R_1)^{\dagger}$,
  where $B^{\dagger}$ denotes the Moore-Penrose pseudoinverse of matrix
  $B$ \cite[pp. 257--258]{Golub1996}. 
\end{algorithmic}
\end{minipage}
\end{flushleft}
\end{small}
\caption{Construct a low-rank approximation $A \approx U_1 M U_2^{\TT}$
  for rank $r$}
\label{alg:randomSVD}
\end{algorithm}

The accuracy of this algorithm and its variants has been studied
thoroughly by several groups. If the matrix $2$-norm is used to
measure the error, it is well-known that the best rank-$r$
approximation is provided by the singular value decomposition
(SVD). When the singular values of $A$ decay rapidly, it has been
shown that Algorithm~\ref{alg:randomSVD} results in almost optimal
factorizations with an overwhelming probability
\cite{HalkoMartinssonTropp:09}. As Algorithm~\ref{alg:randomSVD} is to
be used frequently in our algorithm, we analyze briefly its complexity
step by step. The generation of random numbers is quite efficient,
therefore in practice one may ignore the cost of steps $1$ and
$3$. Step $2$ takes $(r+c)$ \textsf{matvec} of matrix $A$ and
$\Or(n(r+c)^2)$ steps for applying the SVD algorithms on an
$n\times(r+c)$ matrix. The cost of step $4$ is the same as the one of
step $2$. Step $5$ involves the computation of $R_2^{\TT}(A R_1)$,
which takes $\Or(n(r + c)^2)$ steps as we have already computed $A
R_1$ in step $2$. Once $R_2^{\TT}(A R_1)$ is ready, the computation of
$M$ takes additional $\Or((r+c)^3)$ steps. Therefore, the total
complexity of Algorithm~\ref{alg:randomSVD} is $\Or(r + c)$
\textsf{matvec}s plus $\Or(n(r+c)^2)$ extra steps.

\subsection{Top-down construction of $\mc{H}$-matrix}\label{sec:1Dalg}

We illustrate the core idea of our algorithm using a simple
one-dimensional example. The algorithm of constructing a hierarchical
matrix $G$ is a top-down pass. We assume throughout the article that
$G$ is symmetric.

For clarity, we will first consider a one dimension example. The
details of the algorithm in two dimensions will be given in Section 2.
We assume that a symmetric matrix $G$ has a hierarchical low-rank
structure corresponding to a hierarchical dyadic decomposition of the
domain. The matrix $G$ is of dimension $n \times n$ with $n = 2^{L_M}$
for an integer $L_M$. Denote the set for all indices as $\mI_{0;1}$,
where the former subscript indicates the level and the latter is the
index for blocks in each level. At the first level, the set is
partitioned into $\mI_{1;1}$ and $\mI_{1;2}$, with the assumption that
$G(\mI_{1;1},\mI_{1;2})$ and $G(\mI_{1;2},\mI_{1;1})$ are numerically
low-rank, say of rank $r$ for a prescribed error tolerance $\veps$. At
level $l$, each block $\mI_{l-1;i}$ on the above level is dyadically
decomposed into two blocks $\mI_{l;2i-1}$ and $\mI_{l;2i}$ with the
assumption that $G(\mI_{l;2i-1}, \mI_{l;2i})$ and $G(\mI_{l;2i},
\mI_{l;2i-1})$ are also numerically low-rank (with the same rank $r$
for the tolerance $\veps$).  Clearly, at level $l$, we have in total
$2^l$ off-diagonal low-rank blocks.  We stop at level $L_M$, for which
the block $\mI_{L_M,i}$ only has one index $\{i\}$.  For simplicity of
notation, we will abbreviate $G(\mI_{l;i}, \mI_{l;j})$ by $G_{l;ij}$.
We remark that the assumption that off-diagonal blocks are low-rank
matrices may not hold for general elliptic operators in higher
dimensions.  However, this assumption simplifies the introduction of
the concept of our algorithm. More realistic case will be discussed in
detail in Sections~\ref{sec:peeloff} and \ref{sec:peeloffdetails}.

The overarching strategy of our approach is to peel off the
off-diagonal blocks level by level and simultaneously construct their
low-rank approximations. On the first level, $G_{1;12}$ is numerically
low-rank. In order to use the randomized SVD algorithm for $G_{1;12}$,
we need to know the product of $G_{1;12}$ and also $G_{1;12}^{\TT} =
G_{1;21}$ with a collection of random vectors. This can be done by
observing that
\begin{equation}
\begin{pmatrix}
G_{1;11} & G_{1;12} \\
G_{1;21} & G_{1;22} 
\end{pmatrix}
\begin{pmatrix}
R_{1;1} \\
0
\end{pmatrix}
= 
\begin{pmatrix}
G_{1;11}R_{1;1} \\
G_{1;21}R_{1;1}
\end{pmatrix},
\label{eqn:G1R11}
\end{equation}
\begin{equation}
\begin{pmatrix}
G_{1;11} & G_{1;12} \\
G_{1;21} & G_{1;22} 
\end{pmatrix}
\begin{pmatrix}
0 \\
R_{1;2}
\end{pmatrix}
= 
\begin{pmatrix}
G_{1;12}R_{1;2} \\
G_{1;22}R_{1;2}
\end{pmatrix},
\label{eqn:G1R12}
\end{equation}
where $R_{1;1}$ and $R_{1;2}$ are random matrices of dimension
$n/2\times (r+c)$. We obtain $(G_{1;21}R_{1;1})^{\TT} = R_{1;1}^{\TT}
G_{1;12}$ by restricting the right hand side of Eq.~\eqref{eqn:G1R11}
to $\mI_{1;2}$ and obtain $G_{1;12}R_{1;2}$ by restricting the right
hand side of Eq.~\eqref{eqn:G1R12} to $\mI_{1;1}$, respectively. The
low-rank approximation using Algorithm~\ref{alg:randomSVD} results in
\begin{equation}
  G_{1;12} \approx \wh{G}_{1;12} = U_{1;12} M_{1;12} U_{1;21}^\TT.
\end{equation}
$U_{1;12}$ and $U_{1;21}$ are $n/2 \times r$ matrices and
$M_{1;12}$ is an $r\times r$ matrix. Due to the fact that $G$ is
symmetric, a low-rank approximation of $G_{1;21}$ is obtained as the
transpose of $G_{1;12}$.

Now on the second level, the matrix $G$ has the form 
\begin{equation*}
\begin{pmatrix}
\begin{matrix}
G_{2;11} & G_{2;12} \\
G_{2;21} & G_{2;22}
\end{matrix}
& G_{1;12} \\
G_{1;21} &
\begin{matrix}
G_{2;33} & G_{2;34} \\
G_{2;43} & G_{2;44}
\end{matrix}
\end{pmatrix}.
\end{equation*}
The submatrices $G_{2;12}$, $G_{2;21}$, $G_{2;34}$, and $G_{2;43}$
are numerically low-rank, to obtain their low-rank approximations by
the randomized SVD algorithm. Similar to the first level, we could
apply $G$ on random matrices of the form like $(R_{2;1}, 0, 0,
0)^\TT$. This will require $4 (r+c)$ number of matrix-vector
multiplications. However, this is not optimal: Since we already know
the interaction between $\mI_{1;1}$ and $\mI_{1;2}$, we could combine
the calculations together to reduce the number of matrix-vector
multiplications needed.  Observe that
\begin{equation}
  \begin{pmatrix}
    \begin{matrix}
      G_{2;11} & G_{2;12} \\
      G_{2;21} & G_{2;22}
    \end{matrix}
    & G_{1;12} \\
    G_{1;21} &
    \begin{matrix}
      G_{2;33} & G_{2;34} \\
      G_{2;43} & G_{2;44}
    \end{matrix}
  \end{pmatrix}
  \begin{pmatrix}
    R_{2;1} \\ 0 \\ R_{2;3} \\ 0
  \end{pmatrix}
  = 
  \begin{pmatrix}
    \begin{pmatrix}
      G_{2;11} R_{2;1} \\ G_{2;21} R_{2;1}
    \end{pmatrix}
    + G_{1;12}
    \begin{pmatrix}
      R_{2;3} \\ 0
    \end{pmatrix} \\
    \begin{pmatrix}
      G_{2;33} R_{2;3} \\ G_{2;43} R_{2;3}
    \end{pmatrix}
    + G_{1;21}
    \begin{pmatrix}
      R_{2;1} \\ 0
    \end{pmatrix}
  \end{pmatrix}.
\end{equation}
Denote
\begin{equation}
  \wh{G}^{(1)} = 
  \begin{pmatrix}
    0 & \wh{G}_{1;12} \\
    \wh{G}_{1;21} & 0
  \end{pmatrix}
\end{equation}
with $\wh{G}_{1;12}$ and $\wh{G}_{1;21}$ the low-rank approximations
we constructed on the first level, then
\begin{equation}
  \wh{G}^{(1)}\begin{pmatrix}
    R_{2;1} \\ 0 \\ R_{2;3} \\ 0
  \end{pmatrix}
  = \begin{pmatrix}
    \wh{G}_{1;12}
    \begin{pmatrix}
      R_{2;3} \\ 0
    \end{pmatrix} \\
    \wh{G}_{1;21}
    \begin{pmatrix}
      R_{2;1} \\ 0
    \end{pmatrix}
  \end{pmatrix}.
  \label{eqn:combine1d1}
\end{equation}
Therefore,
\begin{equation}
  (G - \wh{G}^{(1)}) \begin{pmatrix}
    R_{2;1} \\ 0 \\ R_{2;3} \\ 0
\end{pmatrix}
\approx
\begin{pmatrix}
  G_{2;11} R_{2;1} \\ G_{2;21} R_{2;1} \\ G_{2;33}R_{2;3} \\
  G_{2;43}R_{2;3}
\end{pmatrix},
\label{eqn:combine1d2}
\end{equation}
so that we simultaneously obtain $(G_{2;21}R_{2;1})^{\TT} =
R_{2;1}^\TT G_{2;12}$ and $(G_{2;43}R_{2;3})^{\TT} = R_{2;3}^\TT G_{2;34}$.
Similarly, applying $G$ on $(0, R_{2;2}, 0, R_{2;4})^\TT$ provides
$G_{2;12}R_{2;2}$ and $G_{2;34}R_{2;4}$.  We can then obtain the
following low-rank approximations by invoking
Algorithm~\ref{alg:randomSVD}.
\begin{equation}
\begin{split}
  G_{2;12} & \approx \wh{G}_{2;12} = U_{2;12} M_{2;12} U_{2;21}^\TT,\\
  G_{2;34} & \approx \wh{G}_{2;34} = U_{2;34} M_{2;34} U_{2;43}^\TT.
\end{split}
\end{equation}
The low-rank approximations of $G_{2;21}$ and $G_{2;43}$ are again
given by the transposes of the above formulas.

Similarly, on the third level, the matrix $G$ has the form 
\begin{equation}
  \begin{pmatrix}
    \begin{matrix}
      \begin{matrix}
        G_{3;11} & G_{3;12} \\
        G_{3;21} & G_{3;22}
      \end{matrix}
      &
      G_{2;12} \\
      G_{2;21} &
      \begin{matrix}
        G_{3;33} & G_{3;34} \\
        G_{3;43} & G_{3;44}
      \end{matrix}
    \end{matrix}&
    G_{1;12} \\
    G_{1;21} &
    \begin{matrix}
      \begin{matrix}
        G_{3;55} & G_{3;56} \\
        G_{3;65} & G_{3;66}
      \end{matrix} &
      G_{2;34} \\
      G_{2;43} &
      \begin{matrix}
        G_{3;77} & G_{3;78} \\
        G_{3;87} & G_{3;88}
      \end{matrix}
    \end{matrix}
  \end{pmatrix},
\end{equation}
and define 
\begin{equation}
  \wh{G}^{(2)} = 
  \begin{pmatrix}
    \begin{matrix}
      0 & \wh{G}_{2;12} \\
      \wh{G}_{2;21} & 0
    \end{matrix} &
    0 \\
    0 &
    \begin{matrix}
      0 & \wh{G}_{2;34} \\
      \wh{G}_{2;43} & 0
    \end{matrix}
  \end{pmatrix}.
\end{equation}
We could simultaneously obtain the product of
$G_{3;12}$, $G_{3;34}$, $G_{3;56}$ and $G_{3;78}$ with random vectors
by applying the matrix $G$ with random vectors of the form 
\begin{equation*}
  (R_{3;1}^\TT,0, R_{3;3}^\TT, 0, R_{3;5}^\TT, 0, R_{3;7}^\TT, 0)^\TT,
\end{equation*}
then subtract the product of $\wh{G}^{(1)} + \wh{G}^{(2)}$ with the
same vectors. Again invoking Algorithm~\ref{alg:randomSVD} provides
us the low-rank approximations of these off-diagonal blocks.

The algorithm continues in the same fashion for higher levels.  The
combined random tests lead to a constant number of \textsf{matvec} at
each level. As there are $\log(n)$ levels in total, the total number
of matrix-vector multiplications scales logarithmically.

When the block size on a level becomes smaller than the given criteria
(for example, the numerical rank $r$ used in the construction), one
could switch to a deterministic way to get the off-diagonal blocks.
In particular, we stop at a level $L$ ($L<L_M$) such that each
$\mI_{L;i}$ contains about $r$ entries. Now only the elements in the
diagonal blocks $G_{L,ii}$ need to be determined. This can be
completed by applying $G$ to the matrix
\[
(I,I,\ldots,I)^{\TT},
\]
where $I$ is the identity matrix whose dimension is equal to the
number of indices in $\mI_{L;i}$.


Let us summarize the structure of our algorithm. From the top level to
the bottom level, we peel off the numerically low-rank off-diagonal
blocks using the randomized SVD algorithm. The matrix-vector
multiplications required by the randomized SVD algorithms are computed
effectively by {\em combining} several random tests into one using the
zero pattern of the {\em remaining} matrix. In this way, we get an
efficient algorithm for constructing the hierarchical representation
for the matrix $G$.

\subsection{Related works}

Our algorithm is built on top of the framework of the
$\mc{H}$-matrices proposed by Hackbusch and his collaborators
\cite{BormGrasedyckHackbusch:06, Hackbusch:99,
  BebendorfHackbusch:03}. The definitions of the $\mc{H}$-matrices
will be summarized in Section~\ref{sec:algorithm}. In a nutshell, the
$\mc{H}$-matrix framework is an operational matrix algebra for
efficiently representing, applying, and manipulating discretizations
of operators from elliptic partial differential equations. Though we
have known how to represent and apply these matrices for quite some
time \cite{GreengardRokhlin:87}, it is the contribution of the
$\mc{H}$-matrix framework that enables one to manipulate them in a
general and coherent way. A closely related matrix algebra is also
developed in a more numerical-linear-algebraic viewpoint under the
name {\em hierarchical semiseparable matrices} by Chandrasekaran, Gu,
and others \cite{ChandrasekaranGuLyons2005,ChandrasekaranGuPals2006}.
Here, we will follow the notations of the $\mc{H}$-matrices as our
main motivations are from numerical solutions of elliptic PDEs.

A basic assumption of our algorithm is the existence of a fast
matrix-vector multiplication subroutine. The most common case is when
$G$ is the inverse of the stiffness matrix $H$ of a general elliptic
operator. Since $H$ is often sparse, much effort has been devoted to
computing $u=Gf$ by solving the linear system $Hu=f$. Many ingenious
algorithms have been developed for this purpose in the past forty
years. Commonly-seen examples include multifrontal algorithms
\cite{DuffReid:83,George:73}, geometric multigrids
\cite{Brandt:77,BriggsHensonMcCormick2000,Hackbusch:99}, algebraic
multigrids (AMG) \cite{Brandt:85}, domain decompositions methods
\cite{SmithBjorstadGropp1996,ToselliWidlund2005}, wavelet-based fast
algorithms \cite{BeylkinCoifmanRokhlin:91} and preconditioned
conjugate gradient algorithms (PCG) \cite{BenziMeyer:96}, to name a
few. Very recently, both Chandrasekaran \etal{}
\cite{ChandrasekaranGuLiXia:06} and Martinsson \cite{Martinsson:09}
have combined the idea of the multifrontal algorithms with the
$\mc{H}$-matrices to obtain highly efficiently direct solvers for
$Hu=f$. Another common case for which a fast matrix-vector
multiplication subroutine is available comes from the boundary
integral equations where $G$ is often a discretization of a Green's
function restricted to a domain boundary. Fast algorithms developed
for this case include the famous fast multipole method
\cite{GreengardRokhlin:87}, the panel clustering method
\cite{HackbuschNowak:89}, and others.  All these fast algorithms
mentioned above can be used as the ``black box'' algorithm for our
method.

As shown in the previous section, our algorithm relies heavily on the
randomized singular value decomposition algorithm for constructing the
factorizations of the off-diagonal blocks. This topic has been a
highly active research area in the past several years and many
different algorithms have been proposed in the literature. Here, for
our purpose, we have adopted the algorithm described in
\cite{LibertyWoolfeMartinssonRokhlinTygert:07,
  WoolfeLibertyRokhlinTygert:08}. In a related but slightly different
problem, the goal is to find low-rank approximations $A=C U R$ where
$C$ contains a subset of columns of $A$ and $R$ contains a subset of
rows. Papers devoted to this task include
\cite{DrineasKannanMahoney2006a,
  DrineasKannanMahoney2006b,Goreinov:97,MahoneyDrineas:09}.  In our
setting, since we assume no direct access of entries of the matrix $A$
but only its impact through matrix-vector multiplications, the
algorithm proposed by \cite{LibertyWoolfeMartinssonRokhlinTygert:07}
is the most relevant choice. An excellent recent review of this fast
growing field can be found in \cite{HalkoMartinssonTropp:09}.

In a recent paper \cite{Martinsson:08}, Martinsson considered also the
problem of constructing the $\mc{H}$-matrix representation of a
matrix, but he assumed that one can access arbitrary entries of the
matrix besides the fast matrix-vector multiplication subroutine. Under
this extra assumption, he showed that one can construct the $\mc{H}^2$
representation of the matrix with $\Or(1)$ matrix-vector
multiplications and accesses of $\Or(n)$ matrix entries. However, in
many situations including the case of $G$ being the inverse of the
stiffness matrix of an elliptic differential operator, accessing
entries of $G$ is by no means a trivial task. Comparing with
Martinsson's work, our algorithm only assumes the existence of a fast
matrix-vector multiplication subroutine, and hence is more general.

As we mentioned earlier, one motivation for computing $G$ explicitly
is to further compress the matrix $G$. The most common example in the
literature of numerical analysis is the process of numerical
homogenization or upscaling \cite{EngquistRunborg:02}. Here the matrix
$G$ is often again the inverse of the stiffness matrix $H$ of an
elliptic partial differential operator. When $H$ contains information
from all scales, the standard homogenization techniques
fail. Recently, Owhadi and Zhang \cite{OwhadiZhang:07} proposed an
elegant method that, under the assumption that the Cordes condition is
satisfied, upscales a general $H$ in divergence form using metric
transformation. Computationally, their approach involves $d$ solves of
form $Hu=f$ with $d$ being the dimension of the problem. On the other
hand, if $G$ is computed using our algorithm, one can obtain the
upscaled operator by inverting a low-passed and down-sampled version of
$G$. Complexity-wise, our algorithm is more costly since it requires
$\Or(\log n)$ solves of $Hu=f$. However, since our approach makes no
analytic assumptions about $H$, it is expected to be more general.


\section{Algorithm}
\label{sec:algorithm}

We now present the details of our algorithm in two dimensions. In
addition to a top-down construction using the peeling idea presented
in the introduction, the complexity will be further reduced using the
$\mc{H}^2$ property of the matrix
\cite{BormGrasedyckHackbusch:06,HackbuschKhoromskijSauter:00}. The
extension to three dimensions is straightforward.

In two dimensions, a more conservative partition of the domain is
required to guarantee the low-rankness of the matrix blocks. We will
start with discussion of this new geometric setup. Then we will recall
the notion of hierarchical matrices and related algorithms in
Section~\ref{sec:Hmatrix}. The algorithm to construct an $\mc{H}^2$
representation for a matrix using matrix-vector multiplications will
be presented in Sections~\ref{sec:peeloff} and
\ref{sec:peeloffdetails}. Finally, variants of the algorithm for
constructing the $\mc{H}^1$ and uniform $\mc{H}^1$ representations
will be described in Section~\ref{sec:peeloffvariants}.

\subsection{Geometric setup and notations}

Let us consider an operator $G$ defined on a 2D domain $[0,1)^2$ with
periodic boundary condition. We discretize the problem using an $n = N
\times N$ uniform grid with $N$ being a power of $2$: $N =
2^{L_M}$. Denote the set of all grid points as
\begin{equation}
  \mI_0 = 
  \{(k_1/N, k_2/N) \mid  k_1, k_2 \in \mathbb{N},\, 0\leq k_1, k_2 < N \}
\end{equation}
and partition the domain hierarchically into $L+1$ levels ($L <
L_M$). On each level $l$ ($0\leq l \leq L$), we have $2^l\times 2^l$
boxes denoted by $\mI_{l;ij} = [(i-1)/2^l, i/2^l)\times[(j-1)/2^l,
j/2^l)$ for $1\le i,j \le 2^l$. The symbol $\mI_{l;ij}$ will also be
used to denote the grid points that lies in the box $\mI_{l;ij}$. The
meaning should be clear from the context. We will also use $\mI_l$(or
$\mJ_l$) to denote a general box on certain level $l$. The subscript
$l$ will be omitted, when the level is clear from the context. For a
given box $\mI_l$ for $l\geq 1$, we call a box $\mJ_{l-1}$ on level
$l-1$ its parent if $\mI_l \subset \mJ_{l-1}$. Naturally, $\mI_l$ is
called a child of $\mJ_{l-1}$. It is clear that each box except those
on level $L$ will have four children boxes.

For any box $\mI$ on level $l$, it covers $N/2^l \times N/2^l$ grid
points. The last level $L$ can be chosen so that the leaf box has a
constant number of points in it (\textit{i.e.} the difference $L_M -
L$ is kept to be a constant when $N$ increases). 

\begin{figure}[ht]
  \begin{center}
    \includegraphics[width=3in]{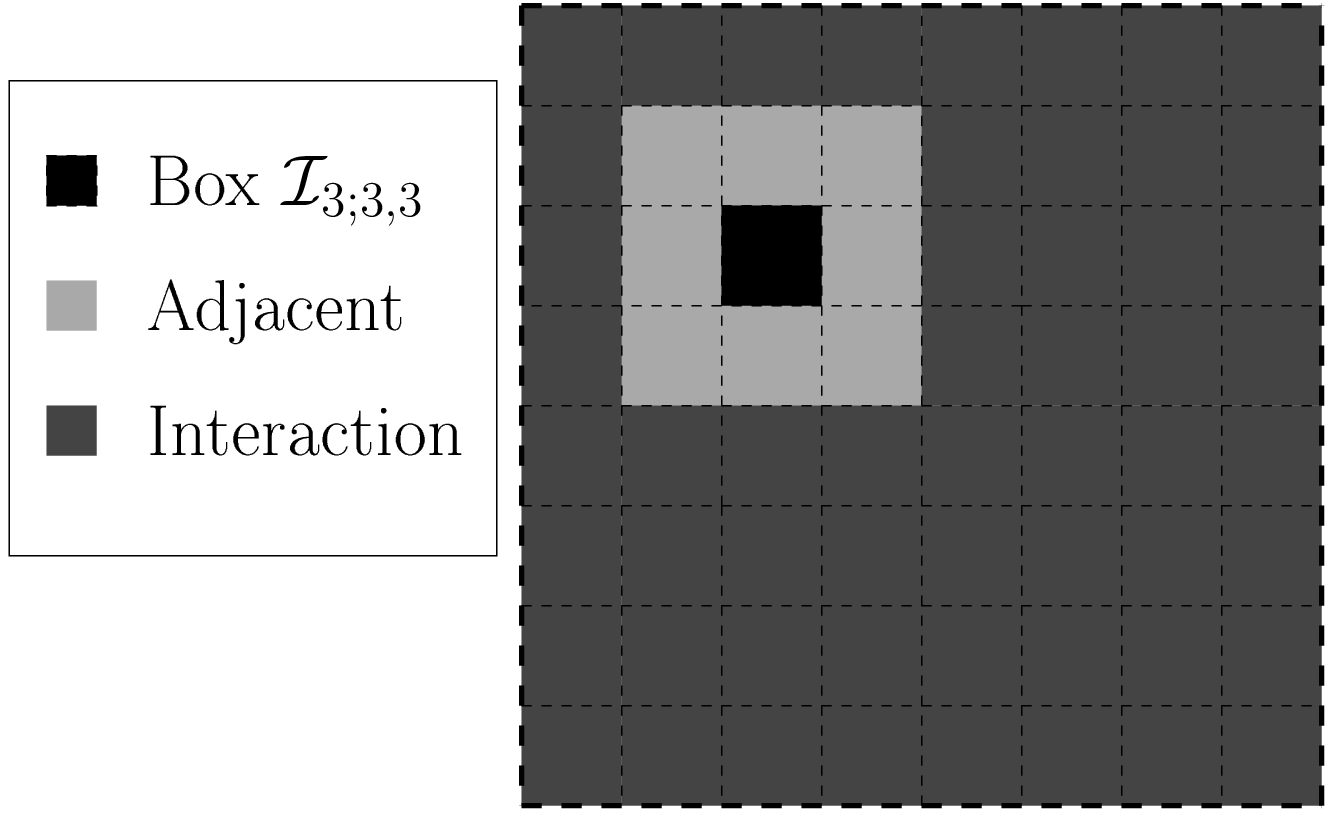}
  \end{center}
  \caption{Illustration of the computational domain at level $3$.
    $\mI_{3;3,3}$ is the black box. The neighbor list
    $\NL(\mI_{3;3,3})$ consists of $8$ adjacent light gray boxes and
    the black box itself, and the interaction list $\IL(\mI_{3;3,3})$
    consists of the $55$ dark gray boxes.}
  \label{fig:level3}
\end{figure}

\begin{figure}[ht]
  \begin{center}
    \includegraphics[width=3in]{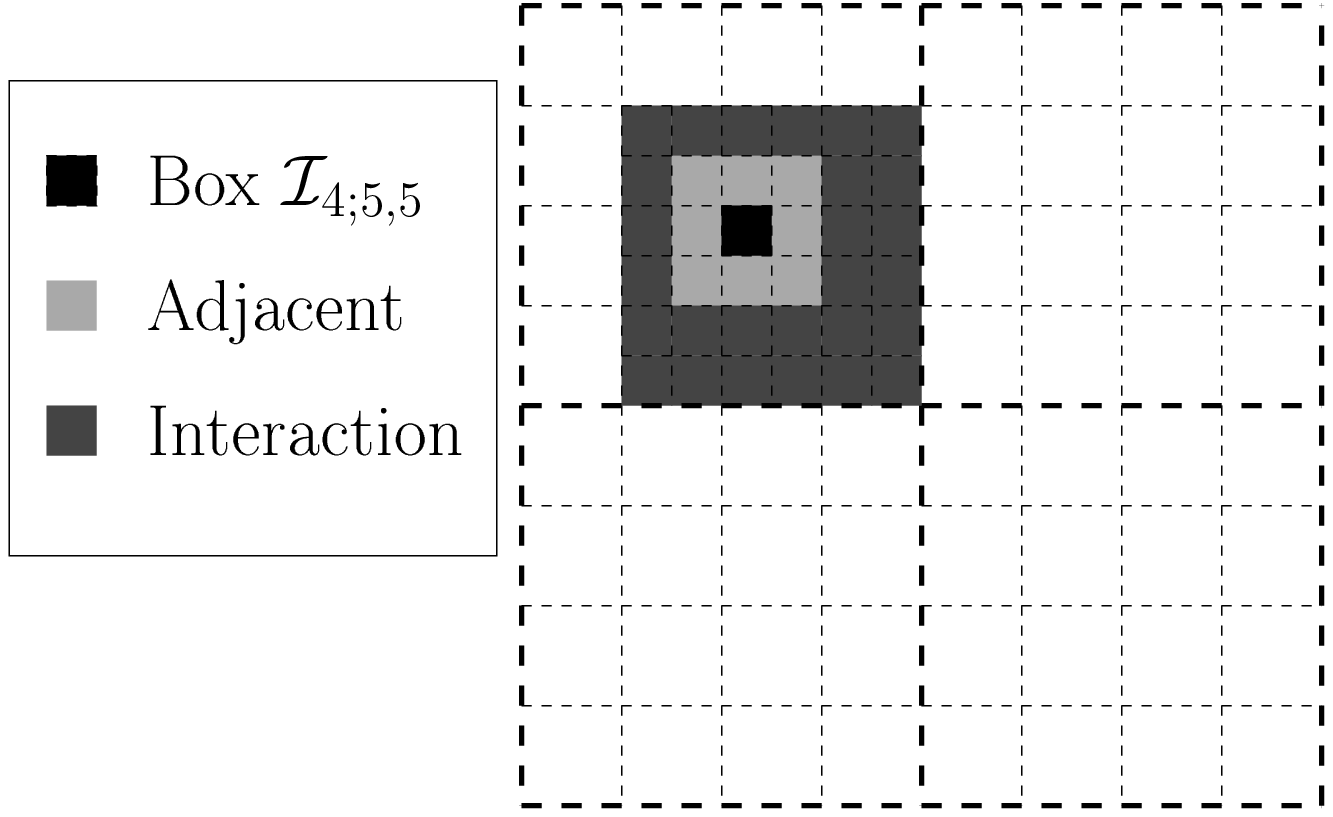}
  \end{center}
  \caption{Illustration of the computational domain at level $4$.
    $\mI_{4;5,5}$ is the black box. The neighbor list
    $\NL(\mI_{4;5,5})$ consists of $8$ adjacent light gray boxes and
    the black box itself, and the interaction list $\IL(\mI_{4;5,5})$
    consists of the $27$ dark gray boxes.}
  \label{fig:level4}
\end{figure}

For simplicity of presentation, we will start the method from level
$3$. It is also possible to start from level $2$. Level $2$ needs to
be treated specially, as for level $3$. We define the following
notations for a box $\mI$ on level $l$ ($l\geq 3$):
\begin{itemize}
\item[$\NL(\mI)$] Neighbor list of box $\mI$. This list contains the
  boxes on level $l$ that are adjacent to $\mI$ and also $\mI$
  itself. There are $9$ boxes in the list for each $\mI$.
\item[$\IL(\mI)$] Interaction list of box $\mI$. When $l = 3$, this
  list contains all the boxes on level $3$ minus the set of boxes in
  $\NL(\mI)$. There are $55$ boxes in total. When $l > 3$, this list
  contains all the boxes on level $l$ that are children of boxes in
  $\NL(\mP)$ with $\mP$ being $\mI$'s parent minus the set of boxes in
  $\NL(\mI)$. There are $27$ such boxes.
\end{itemize}
Notice that these two lists determine two symmetric relationship:
$\mJ\in\NL(\mI)$ if and only if $\mI\in\NL(\mJ)$ and $\mJ\in\IL(\mI)$
if and only if $\mI\in\IL(\mJ)$. Figs.~\ref{fig:level3} and
\ref{fig:level4} illustrate the computational domain and the lists for
$l=3$ and $l=4$, respectively.



\FloatBarrier

For a vector $f$ defined on the $N\times N$ grid $\mI_0$, we define
$f(\mI)$ to be the restriction of $f$ to grid points $\mI$. For a
matrix $G \in \RR^{N^2 \times N^2}$ that represents a linear map from
$\mI_0$ to itself, we define $G(\mI,\mJ)$ to be the restriction of $G$
on $\mI \times \mJ$. 

A matrix $G \in \RR^{N^2 \times N^2}$ has the following decomposition
\begin{equation}\label{eq:Gdecomp}
  G = G^{(3)} + G^{(4)} + \cdots + G^{(L)} + D^{(L)}.
\end{equation}
Here, for each $l$, $G^{(l)}$ incorporates the interaction on level
$l$ between a box with its interaction list. More precisely, $G^{(l)}$
has a $2^{2l} \times 2^{2l}$ block structure:
\begin{equation*}
  G^{(l)}(\mI,\mJ) = 
  \begin{cases}
    G(\mI,\mJ), & \mI \in \IL(\mJ)\ (\text{eq. } \mJ \in \IL(\mI));\\
    0, & \text{otherwise}
  \end{cases}
\end{equation*}
with $\mI$ and $\mJ$ both on level $l$.  The matrix $D^{(L)}$ includes
the interactions between adjacent boxes at level $L$:
\begin{equation*}
  D^{(L)}(\mI,\mJ) = 
  \begin{cases}
    G(\mI,\mJ), & \mI \in \NL(\mJ)\  (\text{eq. } \mJ \in \NL(\mI));\\
    0, & \text{otherwise}
  \end{cases}
\end{equation*}
with $\mI$ and $\mJ$ both on level $L$.   To show that
  \eqref{eq:Gdecomp} is true, it suffices to prove that for any two
  boxes $\mI$ and $\mJ$ on level $L$, the right hand side gives
  $G(\mI, \mJ)$. In the case that $\mI \in \NL(\mJ)$, this is
  obvious. Otherwise, it is clear that we can find a level $l$, and
  boxes $\mI'$ and $\mJ'$ on level $l$, such that $\mI' \in
  \IL(\mJ')$, $\mI \subset \mI'$ and $\mJ \subset \mJ'$, and hence
  $G(\mI, \mJ)$ is given through $G(\mI', \mJ')$. Throughout the
text, we will use $\norm{A}_2$ to denote the matrix $2$-norm of matrix
$A$.


\subsection{Hierarchical matrix}\label{sec:Hmatrix}

Our algorithm works with the so-called hierarchical matrices.  We
recall in this subsection some basic properties of this type of
matrices and also some related algorithms. For simplicity of notations
and representation, we will only work with symmetric matrices. For a
more detailed introduction of the hierarchical matrices and their
applications in fast algorithms, we refer the readers to
\cite{Hackbusch:99, HackbuschKhoromskijSauter:00}.

\subsubsection{$\mc{H}^1$ matrices}
\begin{defn}
  $G$ is a (symmetric) $\mc{H}^1$-matrix if for any $\veps >0$, there
  exists $r(\veps) \lesssim \log(\veps^{-1})$ such that for any pair
  $(\mI, \mJ)$ with $\mI \in \IL(\mJ)$, there exist orthogonal
  matrices $U_{\mI\mJ}$ and $U_{\mJ\mI}$ with $r(\veps)$ columns and
  matrix $M_{\mI\mJ}\in \RR^{r(\veps)\times r(\veps)}$ such that
  \begin{equation}
    \norm{G(\mI,\mJ) - U_{\mI\mJ} M_{\mI\mJ} U_{\mJ\mI}^{\TT} }_2 \leq \veps 
    \norm{G(\mI,\mJ)}_2.
  \label{eqn:H1mat}
  \end{equation}
\end{defn}

The main advantage of the $\mc{H}^1$ matrix is that the application of
such matrix on a vector can be efficiently evaluated: Within error
$\Or(\veps)$, one can use $\wh{G}(\mI,\mJ) = U_{\mI\mJ} M_{\mI\mJ}
U_{\mJ\mI}^{\TT}$, which is low-rank, instead of the original block
$G(\mI,\mJ)$. The algorithm is described in
Algorithm~\ref{alg:H1matvec}. It is standard that the complexity of
the matrix-vector multiplication for an $\mc{H}^1$ matrix is
$\Or(N^2\log N)$ \cite{Hackbusch:99}.

\begin{algorithm}[ht]
\begin{small}
\begin{center}
\begin{minipage}{4.5in}
\begin{algorithmic}[1]
  \STATE $u = 0$; 
  \FOR {$l = 3$ to $L$}
  \FOR {$\mI$ on level $l$}
  \FOR {$\mJ \in \IL(\mI)$}
  \STATE $u(\mI) = u(\mI) + U_{\mI\mJ} (M_{\mI\mJ} (U_{\mJ\mI}^{\TT} f(\mJ)))$;
  \ENDFOR
  \ENDFOR
  \ENDFOR
  \FOR {$\mI$ on level $L$}
  \FOR {$\mJ \in \NL(\mI)$}
  \STATE $u(\mI) = u(\mI) + G(\mI,\mJ) f(\mJ)$;
  \ENDFOR
  \ENDFOR
\end{algorithmic}
\end{minipage}
\end{center}
\end{small}
\caption{Application of a $\mc{H}^1$-matrix $G$ on a vector $f$.}
\label{alg:H1matvec}
\end{algorithm}

\subsubsection{Uniform $\mc{H}^1$ matrix}
\begin{defn}
  $G$ is a (symmetric) uniform $\mc{H}^1$-matrix if for any $\veps
  >0$, there exists $r_U(\veps) \lesssim \log(\veps^{-1})$ such that
  for each box $\mI$, there exists an orthogonal matrix $U_{\mI}$ with
  $r_U(\veps)$ columns such that for any pair $(\mI, \mJ)$ with $\mI
  \in \IL(\mJ)$
  \begin{equation}
    \norm{G(\mI,\mJ) - U_{\mI} N_{\mI\mJ} U_{\mJ}^{\TT} }_2 \leq \veps 
    \norm{G(\mI,\mJ)}_2
    \label{eqn:uniformH1mat}
  \end{equation}
  with $N_{\mI\mJ} \in \RR^{r_U(\veps)\times r_U(\veps)}$.
\end{defn}

The application of a uniform $\mc{H}^1$ matrix to a vector is
described in Algorithm~\ref{alg:unifH1matvec}. The complexity of the
algorithm is still $\Or(N^2\log N)$. However, the prefactor is much
better as each $U_{\mI}$ is applied only once. The speedup over
Algorithm~\ref{alg:H1matvec} is roughly $27 r(\veps) / r_U(\veps)$
\cite{Hackbusch:99}.

\begin{algorithm}[ht]
\begin{small}
\begin{flushleft}
\begin{minipage}{4.8in}
\begin{minipage}[t]{2.35in}
\begin{algorithmic}[1]
  \STATE $u = 0$; 
  \FOR {$l = 3$ to $L$}
  \FOR {$\mJ$ on level $l$}
  \STATE $\wt{f}_{\mJ} = U_{\mJ}^{\TT} f(\mJ)$;
  \ENDFOR
  \ENDFOR
  \FOR {$l = 3$ to $L$} 
  \FOR {$\mI$ on level $l$}
  \STATE $\wt{u}_{\mI} = 0$;
  \FOR {$\mJ \in \IL(\mI)$}
  \STATE $\wt{u}_{\mI} = \wt{u}_{\mI} + N_{\mI\mJ} \wt{f}_{\mJ}$;
  \ENDFOR
  \ENDFOR
  \ENDFOR
\end{algorithmic}
\end{minipage}
\hfill
\begin{minipage}[t]{2.35in}
\begin{algorithmic}[1]
  \setcounter{ALC@line}{14}
  \FOR {$l = 3$ to $L$}
  \FOR {$\mI$ on level $l$}
  \STATE $u(\mI) = u(\mI) + U_{\mI} \wt{u}_{\mI}$;
  \ENDFOR
  \ENDFOR
  \FOR {$\mI$ on level $L$}
  \FOR {$\mJ \in \NL(\mI)$}
  \STATE $u(\mI) = u(\mI) + G(\mI,\mJ) f(\mJ)$;
  \ENDFOR 
  \ENDFOR
\end{algorithmic}
\end{minipage}
\end{minipage}
\end{flushleft}
\end{small}
\caption{Application of a uniform $\mc{H}^1$-matrix $G$ on a vector
  $f$}
\label{alg:unifH1matvec}
\end{algorithm}

\subsubsection{$\mc{H}^2$ matrices}
\begin{defn}
  $G$ is an $\mc{H}^2$ matrix if 
  \begin{itemize}
  \item it is a uniform $\mc{H}^1$ matrix;
  \item Suppose that $\mC$ is any child of a box $\mI$, then
    \begin{equation}
      \norm{U_{\mI}(\mC, :) - U_{\mC} T_{\mC\mI}}_2 \lesssim \veps,
      \label{eqn:H2mat}
    \end{equation}
    for some matrix $T_{\mC\mI} \in \RR^{r_U(\veps)\times
      r_U(\veps)}$.
  \end{itemize}
\end{defn}

The application of an $\mc{H}^2$ matrix to a vector is described in
Algorithm~\ref{alg:H2matvec} and it has a complexity of $\Or(N^2)$,
Notice that, compared with $\mc{H}^1$ matrix, the logarithmic factor
is reduced \cite{HackbuschKhoromskijSauter:00}.

\begin{algorithm}[ht]
\begin{small}
\begin{flushleft}
\begin{minipage}{4.8in}
\begin{minipage}[t]{2.35in}
\begin{algorithmic}[1]
  \STATE $u = 0$; 
  \FOR {$\mJ$ on level $L$}
  \STATE $\wt{f}_{\mJ} = U_{\mJ}^{\TT} f(\mJ)$;
  \ENDFOR
  \FOR {$l = L-1$ down to $3$}
  \FOR {$\mJ$ on level $l$}
  \STATE $\wt{f}_{\mJ} = 0$;
  \FOR {each child $\mC$ of $\mJ$}
  \STATE $\wt{f}_{\mJ} = \wt{f}_{\mJ} + T_{\mC\mJ}^{\TT} \wt{f}_{\mC}$;
  \ENDFOR
  \ENDFOR
  \ENDFOR
  \FOR {$l = 3$ to $L$} 
  \FOR {$\mI$ on level $l$}
  \STATE $\wt{u}_{\mI} = 0$;
  \FOR {$\mJ \in \IL(\mI)$}
  \STATE $\wt{u}_{\mI} = \wt{u}_{\mI} + N_{\mI\mJ} \wt{f}_{\mJ}$;
  \ENDFOR
  \ENDFOR
  \ENDFOR
\end{algorithmic}
\end{minipage}
\hfill
\begin{minipage}[t]{2.35in}
\begin{algorithmic}[1]
  \setcounter{ALC@line}{17}
  \FOR {$l = 3$ to $L-1$}
  \FOR {$\mI$ on level $l$}
  \FOR {each child $\mC$ of $\mI$}
  \STATE $\wt{u}_{\mC} = \wt{u}_{\mC} + T_{\mC\mI}\wt{u}_{\mI}$;
  \ENDFOR
  \ENDFOR
  \ENDFOR
  \FOR {$\mI$ on level $L$}
  \STATE $u(\mI) = U_{\mI} \wt{u}_{\mI}$;
  \ENDFOR
  \FOR {$\mI$ on level $L$}
  \FOR {$\mJ \in \NL(\mI)$}
  \STATE $u(\mI) = u(\mI) + G(\mI,\mJ) f(\mJ)$;
  \ENDFOR 
  \ENDFOR
\end{algorithmic}
\end{minipage}
\end{minipage}
\end{flushleft}
\end{small}
\caption{Application of a $\mc{H}^2$-matrix $G$ on a vector $f$}
\label{alg:H2matvec}
\end{algorithm}

\begin{remark}
  Applying an $\mc{H}^2$ matrix to a vector can indeed be viewed as
  the matrix form of the fast multipole method (FMM)
  \cite{GreengardRokhlin:87}. One recognizes in
  Algorithm~\ref{alg:H2matvec} that the second top-level {\bf for}
  loop corresponds to the M2M (multipole expansion to multipole
  expansion) translations of the FMM; the third top-level {\bf for}
  loop is the M2L (multipole expansion to local expansion)
  translations; and the fourth top-level {\bf for} loop is the L2L
  (local expansion to local expansion) translations.
\end{remark}

In the algorithm to be introduced, we will also need to apply a
partial matrix $G^{(3)} + G^{(4)} + \cdots + G^{(L')}$ for some
$L'\leq L$ to a vector $f$. This amounts to a variant of
Algorithm~\ref{alg:H2matvec}, described in
Algorithm~\ref{alg:partialH2matvec}.

\begin{algorithm}[ht]
\begin{small}
\begin{flushleft}
\begin{minipage}{4.8in}
\begin{minipage}[t]{2.35in}
\begin{algorithmic}[1]
  \STATE $u = 0$; 
  \FOR {$\mJ$ on level $L'$}
  \STATE $\wt{f}_{\mJ} = U_{\mJ}^{\TT} f(\mJ)$;
  \ENDFOR
  \FOR {$l = L'-1$ down to $3$}
  \FOR {$\mJ$ on level $l$}
  \STATE $\wt{f}_{\mJ} = 0$;
  \FOR {each child $\mC$ of $\mJ$}
  \STATE $\wt{f}_{\mJ} = \wt{f}_{\mJ} + T_{\mC\mJ}^{\TT} \wt{f}_{\mC}$;
  \ENDFOR
  \ENDFOR
  \ENDFOR
  \FOR {$l = 3$ to $L'$} 
  \FOR {$\mI$ on level $l$}
  \STATE $\wt{u}_{\mI} = 0$;
  \FOR {$\mJ \in \IL(\mI)$}
  \STATE $\wt{u}_{\mI} = \wt{u}_{\mI} + N_{\mI\mJ} \wt{f}_{\mJ}$;
  \ENDFOR
  \ENDFOR
  \ENDFOR
\end{algorithmic}
\end{minipage}
\hfill
\begin{minipage}[t]{2.35in}
\begin{algorithmic}[1]
  \setcounter{ALC@line}{17}
  \FOR {$l = 3$ to $L'-1$}
  \FOR {$\mI$ on level $l$}
  \FOR {each child $\mC$ of $\mI$}
  \STATE $\wt{u}_{\mC} = \wt{u}_{\mC} + T_{\mC\mI}\wt{u}_{\mI}$;
  \ENDFOR
  \ENDFOR
  \ENDFOR
  \FOR {$\mI$ on level $L'$}
  \STATE $u(\mI) = U_{\mI} \wt{u}_{\mI}$;
  \ENDFOR
\end{algorithmic}
\end{minipage}
\end{minipage}
\end{flushleft}
\end{small}
\caption{Application of a partial $\mc{H}^2$-matrix $G^{(3)} + \cdots
  + G^{(L')}$ on a vector $f$}
\label{alg:partialH2matvec}
\end{algorithm}



\subsection{Peeling algorithm: outline and preparation}\label{sec:peeloff}

We assume that $G$ is a symmetric $\mc{H}^2$ matrix and that there
exists a fast matrix-vector subroutine for applying $G$ to any vector
$f$ as a ``black box''. The goal is to construct an $\mc{H}^2$
representation of the matrix $G$ using only a small
number of test vectors.

The basic strategy is a top-down construction: For each level $l = 3,
\ldots, L$, assume that an $\mc{H}^2$ representation for $G^{(3)} +
\cdots + G^{(l-1)}$ is given, we construct $G^{(l)}$ by the following
three steps:
\begin{enumerate}
\item {\em Peeling}. Construct an $\mc{H}^1$ representation for
  $G^{(l)}$ using the peeling idea and the $\mc{H}^2$ representation
  for $G^{(3)} + \cdots + G^{(l-1)}$.
\item {\em Uniformization.} Construct a uniform $\mc{H}^1$
  representation for $G^{(l)}$ from its $\mc{H}^1$ representation.
\item {\em Projection.} Construct an $\mc{H}^2$ representation for
  $G^{(3)} + \cdots + G^{(l)}$.
\end{enumerate}
The names of these steps will be made clear in the following
discussion.  Variants of the algorithm that only construct an
$\mc{H}^1$ representation (a uniform $\mc{H}^1$ representation,
respectively) of the matrix $G$ can be obtained by only doing the
peeling step (the peeling and uniformization steps,
respectively). These variants will be discussed in Section
\ref{sec:peeloffvariants}.

After we have the $\mc{H}^2$ representation for $G^{(3)} + \cdots +
G^{(L)}$, we use the peeling idea again to extract the diagonal part
$D^{(L)}$. We call this whole process the {\em peeling} algorithm.

Before detailing the peeling algorithm, we mention two procedures that
serve as essential components of our algorithm. The first procedure
concerns with the uniformization step, in which one needs to get a
uniform $\mc{H}^1$ representation for $G^{(l)}$ from its $\mc{H}^1$
representation, \textit{i.e.}, from $\wh{G}(\mI,\mJ) = U_{\mI\mJ}
M_{\mI\mJ} U_{\mJ\mI}^{\TT}$ to $\wh{G}(\mI,\mJ) = U_{\mI} N_{\mI\mJ}
U_{\mJ}^{\TT}$, for all pairs of boxes $(\mI,\mJ)$ with $\mI \in
\IL(\mJ)$. To this end, what we need to do is to find the column space
of
\begin{equation}
  [U_{\mI\mJ_1}M_{\mI\mJ_1} \mid U_{\mI\mJ_2}M_{\mI\mJ_2} 
  \mid \cdots \mid U_{\mI\mJ_t}M_{\mI\mJ_t}],
  \label{eqn:columnspace}
\end{equation}
where $\mJ_j$ are the boxes in $\IL(\mI)$ and $t =
\abs{\IL(\mI)}$. Notice that we weight the singular vectors $U$ by
$M$, so that the singular vectors corresponding to larger singular
values will be more significant. This column space can be found by the
usual SVD algorithm or a more effective randomized version presented
in Algorithm~\ref{alg:H1tounifH1}.  The important left singular vectors
are denoted by $U_{\mI}$, and the diagonal matrix formed by the singular
values associated with $U_{\mI}$ is denoted by $S_{\mI}$.
\begin{algorithm}[ht]
\begin{small}
\begin{flushleft}
\begin{minipage}{4.8in}
\begin{algorithmic}[1]
  \FOR {each box $\mI$ on level $l$}
  \STATE Generate a Gaussian random matrix $R \in \RR^{(r(\veps)\times
    t)\times (r_U(\veps) + c)}$;
  \STATE Form product $[U_{\mI\mJ_1}M_{\mI\mJ_1} \mid 
  \cdots \mid U_{\mI\mJ_t}M_{\mI\mJ_t}]R$ and apply SVD to it. 
  The first $r_U(\veps)$ left singular vectors give $U_{\mI}$, and the
  corresponding singular values give a diagonal matrix $S_{\mI}$;
  \FOR {$\mJ_j \in \IL(\mI)$}
  \STATE $I_{\mI\mJ_j} = U_{\mI}^{\TT} U_{\mI\mJ_j}$;
  \ENDFOR
  \ENDFOR
  \FOR {each pair $(\mI, \mJ)$ on level $l$ with $\mI \in \IL(\mJ)$}
  \STATE $N_{\mI\mJ} = I_{\mI\mJ} M_{\mI\mJ} I_{\mJ\mI}^{\TT}$;
  \ENDFOR
\end{algorithmic}
\end{minipage}
\end{flushleft}
\end{small}
\caption{Construct a uniform $\mc{H}^1$ representation of $G$ from the
  $\mc{H}^1$ representation at a level $l$}
\label{alg:H1tounifH1}
\end{algorithm}

Complexity analysis: For a box $\mI$ on level $l$, the number of grid
points in $\mI$ is $(N/2^l)^2$. Therefore, $U_{\mI\mJ_j}$ are all of
size $(N/2^l)^2\times r(\veps)$ and $M_{\mI\mJ}$ are of size
$r(\veps)\times r(\veps)$. Forming the product
$[U_{\mI\mJ_1}M_{\mI\mJ_1} \mid \cdots \mid
U_{\mI\mJ_t}M_{\mI\mJ_t}]R$ takes $\Or( (N/2^l)^2 r(\veps) (r_U(\veps)
+ c))$ steps and SVD takes $\Or( (N/2^l)^2 (r_U(\veps) + c)^2 )$
steps.  As there are $2^{2l}$ boxes on level $l$, the overall cost of
Algorithm~\ref{alg:H1tounifH1} is $\Or(N^2 (r_U(\veps) + c)^2) =
\Or(N^2)$. One may also apply to $[U_{\mI\mJ_1}M_{\mI\mJ_1} \mid
\cdots \mid U_{\mI\mJ_t}M_{\mI\mJ_t}]$ the deterministic SVD
algorithm, which has the same order of complexity but with a prefactor
about $27 r(\veps) / (r_U(\veps) + c)$ times larger.

The second procedure is concerned with the projection step of the
above list, in which one constructs an $\mc{H}^2$ representation for
$G^{(3)} + \cdots G^{(l)}$. Here, we are given the $\mc{H}^2$
representation for $G^{(3)} + \cdots + G^{(l-1)}$ along with the
uniform $\mc{H}^1$ representation for $G^{(l)}$ and the goal is to
compute the transfer matrix $T_{\mC\mI}$ for a box $\mI$ on level
$l-1$ and its child $\mC$ on level $l$ such that
\[
\norm{U_{\mI}(\mC, :) - U_{\mC} T_{\mC\mI}}_2 \lesssim \veps.
\]
In fact, the existing $U_{\mC}$ matrix of the uniform $\mc{H}^1$
representation may not be rich enough to contain the columns of
$U_{\mI}(\mC, :)$ in its span. Therefore, one needs to update the
content of $U_{\mC}$ as well. To do that, we perform a singular value
decomposition for the combined matrix
\[
[ U_{\mI}(\mC,:) S_{\mI} \mid U_{\mC} S_{\mC}]
\]
and define a matrix $V_{\mC}$ to contain $r_U(\veps)$ left singular
vectors. Again $U_{\mI}, U_{\mC}$ should be weighted by the corresponding
singular values.
The transfer matrix $T_{\mC\mI}$ is then given by
\[
T_{\mC\mI} = V_{\mC}^{\TT} U_{\mI}(\mC, :)
\]
and the new $U_{\mC}$ is set to be equal to $V_{\mC}$. Since $U_{\mC}$
has been changed, the matrices $N_{\mC\mD}$ for $\mD \in \IL(\mC)$
and also the corresponding singular values $S_{\mC}$ need
  to be updated as well.  The details are listed in
Algorithm~\ref{alg:unifH1toH2}.
\begin{algorithm}[ht]
\begin{small}
\begin{flushleft}
\begin{minipage}{4.8in}
\begin{algorithmic}[1]
  \FOR {each box $\mI$ on level $l-1$} 
  \FOR {each child $\mC$ of $\mI$}
  \STATE Form matrix $[U_{\mI}(\mC,:) S_{\mI} \mid U_{\mC} S_{\mC}]$
  and apply SVD to it. The first $r_U(\veps)$ left singular vectors
  give $V_{\mC}$, and the corresponding singular values give a
  diagonal matrix $W_{\mC}$;
  \STATE $K_{\mC} = V_{\mC}^{\TT} U_{\mC}$;
  \STATE $T_{\mC\mI} = V_{\mC}^{\TT} U_{\mI}(\mC, :)$;
  \STATE $U_{\mC} = V_{\mC}$;
  \STATE $S_{\mC} = W_{\mC}$;
  \ENDFOR
  \ENDFOR

  \FOR {each pair $(\mC, \mD)$ on level $l$ with $\mC \in \IL(\mD)$}
  \STATE $N_{\mC\mD} = K_{\mC} N_{\mC\mD} K_{\mD}^{\TT}$;
  \ENDFOR
\end{algorithmic}
\end{minipage}
\end{flushleft}
\end{small}
\caption{Construct an $\mc{H}^2$ representation of $G$ from the
  uniform $\mc{H}^1$ representation at level $l$}
\label{alg:unifH1toH2}
\end{algorithm}

Complexity analysis: The main computational task of
Algorithm~\ref{alg:unifH1toH2} is again the SVD part. For a box $\mC$
on level $l$, the number of grid points in $\mI$ is $(N/2^l)^2$.
Therefore, the combined matrix $[U_{\mI}(\mC,:) S_{\mI} \mid U_{\mC} S_{\mC}]$ is of
size $(N/2^l)^2 \times 2r_U(\veps)$. The SVD computation clearly takes
$\Or((N/2^l)^2 r_U(\veps)^2) = \Or((N/2^l)^2)$ steps. Taking into the
consideration that there are $2^{2l}$ boxes on level $l$ gives rise to
an $\Or(N^2)$ estimate for the cost of Algorithm~\ref{alg:unifH1toH2}.

\subsection{Peeling algorithm: details}\label{sec:peeloffdetails}

With the above preparation, we are now ready to describe the peeling
algorithm in detail at different levels, starting from level $3$.  At
each level, we follow exactly the three steps listed at the beginning
of Section~\ref{sec:peeloff}.

\subsubsection{Level $3$}

First in the peeling step, we construct the $\mc{H}^1$ representation
for $G^{(3)}$. For each pair $(\mI, \mJ)$ on level $3$ such that $\mI
\in \IL(\mJ)$, we will invoke randomized SVD
Algorithm~\ref{alg:randomSVD} to construct the low rank approximation
of $G_{\mI,\mJ}$. However, in order to apply the algorithm we need to
compute $G(\mI,\mJ) R_{\mJ}$ and $R_{\mI}^{\TT} G(\mI,\mJ)$, where
$R_{\mI}$ and $R_{\mJ}$ are random matrices with $r(\veps) + c$
columns. For each box $\mJ$ on level $3$, we construct a matrix $R$ of
size $N^2 \times (r(\veps) + c)$ such that
\begin{equation*}
  R(\mJ, :) = R_{\mJ} \quad \text{and}\quad R(\mI_0\backslash\mJ,:) = 0.
\end{equation*}
Computing $GR$ using $r(\veps) + c$ \textsf{matvec}s and restricting
the result to grid points $\mI \in \IL(\mJ)$ gives $G(\mI,\mJ)
R_{\mJ}$ for each $\mI \in \IL(\mJ)$. 

After repeating these steps for all boxes on level $3$, we hold for
any pair $(\mI, \mJ)$ with $\mI \in \IL(\mJ)$ the following data:
\begin{equation*}
  G(\mI,\mJ) R_{\mJ} \quad\text{and}\quad R_{\mI}^{\TT} G(\mI,\mJ) = (G(\mJ,\mI)
  R_{\mI})^{\TT}.
\end{equation*}
Now, applying Algorithm~\ref{alg:randomSVD} to them gives the low-rank
approximation
\begin{equation}
  \wh{G}(\mI,\mJ) = U_{\mI\mJ} M_{\mI\mJ} U_{\mJ\mI}^{\TT}.
\end{equation}

In the uniformization step, in order to get the uniform $\mc{H}^1$
representation for $G^{(3)}$, we simply apply
Algorithm~\ref{alg:H1tounifH1} to the boxes on level $3$ to get the
approximations
\begin{equation}
  \wh{G}(\mI,\mJ) = U_{\mI} N_{\mI\mJ} U_{\mJ}^{\TT}.
\end{equation}

Finally in the projection step, since we only have $1$ level now
(level $3$), we have already the $\mc{H}^2$ representation for
$G^{(3)}$.

Complexity analysis: The dominant computation is the construction of
the $\mc{H}^1$ representation for $G^{(3)}$. This requires
$r(\veps)+c$ \textsf{matvec}s for each box $\mI$ on level $3$. Since
there are in total $64$ boxes at this level, the total cost is
$64(r(\veps) + c)$ \textsf{matvec}s. From the complexity analysis in
Section~\ref{sec:peeloff}, the computation for the second and third
steps cost an extra $\Or(N^2)$ steps.

\subsubsection{Level $4$}
First in the peeling step, in order to construct the $\mc{H}^1$
representation for $G^{(4)}$, we need to compute the matrices
$G(\mI,\mJ)R_{\mJ}$ and $R_{\mI}^{\TT} G(\mI,\mJ)$ for each pair
$(\mI, \mJ)$ on level $4$ with $\mI\in \IL(\mJ)$. Here $R_{\mI}$ and
$R_{\mJ}$ are again random matrices with $r(\veps) + c$ columns.

One approach is to follow exactly what we did for level $3$: Fix a box
$\mJ$ at this level, construct $R$ of size $N^2\times (r(\veps)+c)$
such that
\begin{equation*}
  R(\mJ, :) = R_{\mJ} \quad \text{and}\quad R(\mI_0\backslash\mJ,:) = 0.
\end{equation*}
Next apply $G - G^{(3)}$ to $R$, by subtracting $GR$ and
$G^{(3)}R$. The former is computed using $r(\veps) + c$
\textsf{matvec}s and the latter is done by
Algorithm~\ref{alg:partialH2matvec}. Finally, restrict the result to
grid points $\mI \in \IL(\mJ)$ gives $G(\mI,\mJ) R_{\mJ}$ for each
$\mI \in \IL(\mJ)$.

However, we have observed in the simple one-dimensional example in
Section~\ref{sec:1Dalg} that random tests can be combined together as
in Eq.~\eqref{eqn:combine1d1} and \eqref{eqn:combine1d2}. We shall
detail this observation in the more general situation here as
following.  Observe that $G - G^{(3)} = G^{(4)} + D^{(4)}$, and
$G^{(4)}(\mJ, \mI)$ and $D^{(4)}(\mJ, \mI)$ for boxes $\mI$ and $\mJ$
on level $4$ is only nonzero if $\mI \in \NL(\mJ) \cup \IL(\mJ)$.
Therefore, $(G - G^{(3)})R$ for $R$ coming from $\mJ$ can only be
nonzero in $\NL(\mP)$ with $\mP$ being $\mJ$'s parent. The rest is
automatically zero (up to error $\veps$ as $G^{(3)}$ is approximated
by its $\mc{H}^2$ representation).  Therefore, we can combine the
calculation of different boxes as long as their non-zero regions do
not overlap.

\begin{figure}[ht]
  \begin{center}
    \includegraphics[width=3in]{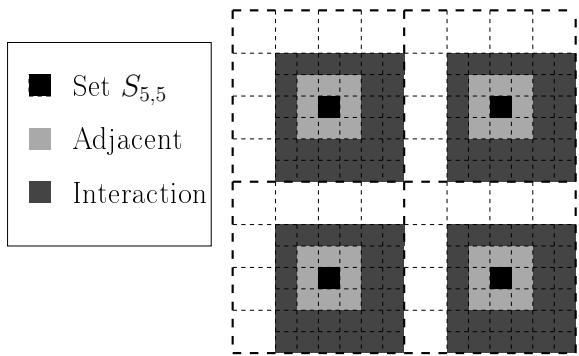}
  \end{center}
  \caption{Illustration of the set $S_{55}$ at level $4$.  This set
    consists of four black boxes
    $\{\mI_{4;5,5},\mI_{4;13,5},\mI_{4;5,13},\mI_{4;13,13}\}$. The
    light gray boxes around each black box are in the neighbor list
    and the dark gray boxes in the interaction list. }
  \label{fig:level4combined}
\end{figure}

More precisely, we introduce the following sets $\mc{S}_{pq}$ for $1\leq
p,q\le 8$ with
\begin{equation}
  \mc{S}_{pq} = \{ \mJ_{4;ij} \mid i \equiv p\,(\text{mod }8),\; j \equiv q\,
  (\text{mod }8) \}.
\end{equation}
There are $64$ sets in total, each consisting of four boxes.
Fig.~\ref{fig:level4combined} illustrates one such set at level $4$.
For each set $\mc{S}_{pq}$, first construct $R$ with
\begin{equation*}
  R(\mJ, :) = 
  \begin{cases}
    R_{\mJ}, & \mJ \in \mc{S}_{pq}; \\
    0, & \text{otherwise}.
  \end{cases}
\end{equation*}
Then, we apply $G - G^{(3)}$ to $R$, by subtracting $GR$ and
$G^{(3)}R$. The former is computed using $r(\veps) + c$
\textsf{matvec}s and the latter is done by
Algorithm~\ref{alg:partialH2matvec}. For each $\mJ\in \mc{S}_{pq}$,
restricting the result to $\mI \in \IL(\mJ)$ gives
$G(\mI,\mJ)R_{\mJ}$. Repeating this computation for all sets $\mc{S}_{pq}$
then provides us with the following data:
\begin{equation*}
  G(\mI,\mJ) R_{\mJ} \quad\text{and}\quad R_{\mI}^{\TT} G(\mI,\mJ) = (G(\mJ,\mI)
  R_{\mI})^{\TT},
\end{equation*}
for each pair $(\mI, \mJ)$ with $\mI \in \IL(\mJ)$. Applying
Algorithm~\ref{alg:randomSVD} to them gives the required low-rank
approximations
\begin{equation}
  \wh{G}(\mI,\mJ) = U_{\mI\mJ} M_{\mI\mJ} U_{\mJ\mI}^{\TT}
\end{equation}
with $U_{\mI\mJ}$ orthogonal.

Next in the uniformization step, the task is to construct the uniform
$\mc{H}^1$ representation of $G^{(4)}$. Similar to the computation
at level $3$, we simply apply Algorithm~\ref{alg:H1tounifH1} to the
boxes on level $4$ to get 
\begin{equation}
  \wh{G}(\mI,\mJ) = U_{\mI} N_{\mI\mJ} U_{\mJ}^{\TT}.
\end{equation}

Finally in the projection step, to get $\mc{H}^2$ representation for
$G^{(3)} + G^{(4)}$, we invoke Algorithm~\ref{alg:unifH1toH2} at level
$4$. Once it is done, we hold the transfer matrices $T_{\mC\mI}$
between any $\mI$ on level $3$ and each of its children $\mC$, along
with the updated uniform $\mc{H}^1$-matrix representation of
$G^{(4)}$.

Complexity analysis: The dominant computation is again the
construction of $\mc{H}^1$ representation for $G^{(4)}$. For each
group $\mc{S}_{pq}$, we apply $G$ to $r(\veps) + c$ vectors and apply
$G^{(3)}$ to $r(\veps) + c$ vectors. The latter takes $\Or(N^2)$ steps
for each application.  Since there are $64$ sets in total, this
computation takes $64(r(\veps) + c)$ \textsf{matvec}s and $\Or(N^2)$
extra steps.

\subsubsection{Level $l$}
First in the peeling step, to construct the $\mc{H}^1$ representation
for $G^{(l)}$, we follow the discussion of level $4$. Define $64$ sets
$\mc{S}_{pq}$ for $1\leq p, q \leq 8$ with
\begin{equation}
  \mc{S}_{pq} = \{ \mJ_{l;ij} \mid i \equiv p\,(\text{mod }8),\; j \equiv q\,
  (\text{mod }8) \}.
\end{equation}
Each set contains exactly $2^l/8 \times 2^l/8$ boxes. For each set
$\mc{S}_{pq}$, construct $R$ with
\begin{equation*}
  R(\mJ, :) = 
  \begin{cases}
    R_{\mJ}, & \mJ \in \mc{S}_{pq}; \\
    0, & \text{otherwise}.
  \end{cases}
\end{equation*}
Next, apply $G - [G^{(3)} + \cdots + G^{(l-1)}]$ to $R$, by
subtracting $GR$ and $[G^{(3)} + \cdots + G^{(l-1)}]R$. The former is
again computed using $r(\veps) + c$ \textsf{matvec}s and the latter is
done by Algorithm~\ref{alg:partialH2matvec} using the $\mc{H}^2$
representation of $G^{(3)} + \cdots + G^{(l-1)}$. For each $\mJ\in
\mc{S}_{pq}$, restricting the result to $\mI \in \IL(\mJ)$ gives
$G(\mI,\mJ)R_{\mJ}$.  Repeating this computation for all sets $\mc{S}_{pq}$
gives the following data for any pair $(\mI, \mJ)$ with $\mI \in
\IL(\mJ)$
\begin{equation*}
  G(\mI,\mJ) R_{\mJ} \quad\text{and}\quad R_{\mI}^{\TT} G(\mI,\mJ) = (G(\mJ,\mI) 
  R_{\mI})^{\TT}.
\end{equation*}
Now applying Algorithm~\ref{alg:randomSVD} to them gives the low-rank
approximation
\begin{equation}
  \wh{G}(\mI,\mJ) = U_{\mI\mJ} M_{\mI\mJ} U_{\mJ\mI}^{\TT}
\end{equation}
with $U_{\mI\mJ}$ orthogonal.

Similar to the computation at level $4$, the uniformization step
that constructs the uniform $\mc{H}^1$ representation of $G^{(l)}$
simply by Algorithm~\ref{alg:H1tounifH1} to the boxes on level
$l$. The result gives the approximation
\begin{equation}
  \wh{G}(\mI,\mJ) = U_{\mI} N_{\mI\mJ} U_{\mJ}^{\TT}.
\end{equation}

Finally in the projection step, one needs to compute an $\mc{H}^2$
representation for $G^{(3)} + \cdots + G^{(l)}$. To this end, we apply
Algorithm~\ref{alg:unifH1toH2} to level $l$.


The complexity analysis at level $l$ follows exactly the one of level
$4$. Since we still have exactly $64$ sets $\mc{S}_{pq}$, the computation
again takes $64(r(\veps) + c)$ \textsf{matvec}s along with $\Or(N^2)$
extra steps. 

These three steps (peeling, uniformization, and projection) are
repeated for each level until we reach level $L$. At this point, we
hold the $\mc{H}^2$ representation for $G^{(3)} + \cdots G^{(L)}$.

\subsubsection{Computation of $D^{(L)}$}
Finally we construct of the diagonal part
\begin{equation}
  D^{(L)} = G - (G^{(3)} + \cdots + G^{(L)}).
\end{equation}
More specifically, for each box $\mJ$ on level $L$, we need to compute
$G(\mI,\mJ)$ for $\mI \in \NL(\mJ)$.

Define a matrix $E$ of size $N^2 \times (N/2^L)^2$ (recall that the
box $\mJ$ on level $L$ covers $(N/2^L)^2$ grid points) by
\begin{equation*}
  E(\mJ, :) = I \quad \text{and} \quad E(\mI_0\backslash\mJ, :) = 0,
\end{equation*}
where $I$ is the $(N/2^L)^2 \times (N/2^L)^2$ identity
matrix. Applying $G - (G^{(3)} + \cdots + G^{(L)})$ to $E$ and
restricting the results to $\mI \in \NL(\mJ)$ gives $G(\mI,\mJ)$ for
$\mI \in \NL(\mJ)$. However, we can do better as $(G - (G^{(3)} +
\cdots + G^{(L)}))E $ is only non-zero in $\NL(\mJ)$. Hence, one can
combine computation of different boxes $\mJ$ as long as $\NL(\mJ)$ do
not overlap.

To do this, define the following $4\times 4 = 16$ sets $\mc{S}_{pq}$,
$1\leq p, q\le 4$
\begin{equation*}
  \mc{S}_{pq} = \{\mJ_{L,ij} \mid  i \equiv p\,(\text{mod }4),\; j \equiv q\,
  (\text{mod }4) \}. 
\end{equation*}
For each set $\mc{S}_{pq}$, construct matrix $E$ by
\begin{equation*}
  E(\mJ, :) = 
  \begin{cases}
    I, & \mJ \in \mc{S}_{pq}; \\
    0, & \text{otherwise}.
  \end{cases}
\end{equation*}
Next, apply $G - (G^{(3)} + \cdots + G^{(L)})$ to $E$. For each
$\mJ\in \mc{S}_{pq}$, restricting the result to $\mI\in \NL(\mJ)$ gives
$G(\mI,\mJ) I = G(\mI,\mJ)$. Repeating this computation for all $16$
sets $\mc{S}_{pq}$ gives the diagonal part $D^{(L)}$.

Complexity analysis: The dominant computation is for each group
$\mc{S}_{pq}$ apply $G$ and $G^{(3)} + \cdots + G^{(L)}$ to $E$, the former
takes $(N/2^L)^2$ \textsf{matvec}s and the latter takes $\Or((N/2^L)^2
N^2)$ extra steps. Recall by the choice of $L$, $N/2^L$ is a
constant. Therefore, the total cost for $16$ sets is $16 (N/2^L)^2 =
\Or(1)$ \textsf{matvec}s and $\Or(N^2)$ extra steps.

Let us now summarize the complexity of the whole peeling
algorithm. From the above discussion, it is clear that at each level
the algorithm spends $64(r(\veps) + c) = \Or(1)$ \textsf{matvec}s and
$\Or(N^2)$ extra steps. As there are $\Or(\log N)$ levels, the overall
cost of the peeling algorithm is equal to $\Or(\log N)$
\textsf{matvec}s plus $\Or(N^2 \log N)$ steps.


 It is a natural concern that whether the error from
  low-rank decompositions on top levels accumulates in the peeling
  steps. As observed from numerical examples in
  Section~\ref{sec:numerical}, it does not seem to be a problem at least for the examples considered. We do not have a
  proof for this though. 

\subsection{Peeling algorithm: variants}\label{sec:peeloffvariants}

In this section, we discuss two variants of the peeling algorithm. Let
us recall that the above algorithm performs the following three steps
at each level $l$.
\begin{enumerate}
\item {\em Peeling}. Construct an $\mc{H}^1$ representation for
  $G^{(l)}$ using the peeling idea and the $\mc{H}^2$ representation
  for $G^{(3)} + \cdots + G^{(l-1)}$.
\item {\em Uniformization.} Construct a uniform $\mc{H}^1$
  representation for $G^{(l)}$ from its $\mc{H}^1$ representation.
\item {\em Projection.} Construct an $\mc{H}^2$ representation for
  $G^{(3)} + \cdots + G^{(l)}$.
\end{enumerate}
As this algorithm constructs the $\mc{H}^2$ representation of the
matrix $G$, we also refer to it more specifically as the $\mc{H}^2$
version of the peeling algorithm. In what follows, we list two simpler
versions that are useful in practice
\begin{itemize}
\item the $\mc{H}^1$ version, and
\item the uniform $\mc{H}^1$ version.
\end{itemize}

In the $\mc{H}^1$ version, we only perform the peeling step at each
level. Since this version constructs only the $\mc{H}^1$
representation, it will use the $\mc{H}^{1}$ representation of
$G^{(3)} + \cdots + G^{(l)}$ in the computation of $(G^{(3)} + \cdots
+ G^{(l)}) R$ within the peeling step at level $l+1$.

In the uniform $\mc{H}^1$ version, we perform the peeling step and the
uniformization step at each level. This will give us instead the
uniform $\mc{H}^1$ version of the matrix. Accordingly, one needs to
use the uniform $\mc{H}^{1}$ representation of $G^{(3)} + \cdots +
G^{(l)}$ in the computation of $(G^{(3)} + \cdots + G^{(l)}) R$
within the peeling step at level $l+1$.


These two simplified versions are of practical value since there are
matrices that are in the $\mc{H}^1$ or the uniform $\mc{H}^1$ class but
not the $\mc{H}^2$ class. A simple calculation shows that these two
simplified versions still take $\Or(\log N)$ \textsf{matvec}s but
requires $\Or(N^2 \log^2 N)$ extra steps. Clearly, the number of extra
steps is $\log N$ times more expensive than the one of the $\mc{H}^2$
version. However, if the fast matrix-vector multiplication subroutine
itself takes $\Or(N^2 \log N)$ steps per application, using the
$\mc{H}^1$ or the uniform $\mc{H}^1$ version does not change the
overall asymptotic complexity.

Between these two simplified versions, the uniform $\mc{H}^1$ version
requires the uniformization step, while the $\mc{H}^1$ version does
not. This seems to suggest that the uniform $\mc{H}^1$ version is more
expensive. However, because (1) our algorithm also utilizes the
partially constructed representations for the calculation at future
levels and (2) the uniform $\mc{H}^1$ representation is much faster to
apply, the construction of the uniform $\mc{H}^{1}$ version turns out
to be much faster. Moreover, since the uniform $\mc{H}^1$
representation stores one $U_{\mI}$ matrix for each box $\mI$ while
the $\mc{H}^1$ version stores about $27$ of them, the uniform
$\mc{H}^1$ is much more memory-efficient, which is very important for
problems in higher dimensions.

\section{Numerical results}\label{sec:numerical}

We study the performance of the hierarchical matrix construction
algorithm for the inverse of a discretized elliptic operator.  The
computational domain is a two dimensional square $[0,1)^2$ with
periodic boundary condition, discretized as an $N\times N$ equispaced
grid.  We first consider the operator $H=-\Delta + V$ with $\Delta$
being the discretized Laplacian operator and the potential being
$V(i,j) = 1 + W(i,j), \quad i,j=1,\ldots,N$.  For all $(i,j)$,
$W(i,j)$ are independent random numbers uniformly distributed in
$[0,1]$.  The potential function $V$ is chosen to have this strong
randomness in order to show that the existence of $\mc{H}$-matrix
representation of the Green's function depends weakly on the
smoothness of the potential.  The inverse matrix of $H$ is denoted by
$G$.  The algorithms are implemented using \MATLAB.  All numerical
tests are carried out on a single-CPU machine.

We analyze the performance statistics by examining both the cost and
the accuracy of our algorithm. The cost factors include the time cost
and the memory cost. While the memory cost is mainly determined by how
the matrix $G$ is compressed and does not depend much on the
particular implementation, the time cost depends heavily on the
performance of \textsf{matvec} subroutine. Therefore, we report both
the wall clock time consumption of the algorithm and the number of
calls to the \textsf{matvec} subroutine. The \textsf{matvec}
subroutine used here is a nested dissection reordered block Gauss
elimination method \cite{George:73}.
For an $N\times N$ discretization of the computational domain, this
\textsf{matvec} subroutine has a computational cost of $\Or(N^2\log
N)$ steps. 

Table~\ref{tab:time} summarizes the \textsf{matvec} number, and the
time cost per degree of freedom (DOF) for the $\mc{H}^1$, the uniform
$\mc{H}^1$ and the $\mc{H}^2$ representations of the peeling algorithm.  The
time cost per DOF is defined by the total time cost divided by the
number of grid points $N^2$.  For the $\mc{H}^1$ and the uniform
$\mc{H}^1$ versions, the error criterion $\veps$ in
Eq.~\eqref{eqn:H1mat}, Eq.~\eqref{eqn:uniformH1mat} and
Eq.~\eqref{eqn:H2mat} are all set to be $10^{-6}$.

The number of calls to the \textsf{matvec} subroutine is the same in
all three cases (as the peeling step is the same for all cases) and is
reported in the third column of Table~\ref{tab:time}.  It is confirmed
that the number of calls to \textsf{matvec} increases logarithmically
with respect to $N$ if the domain size at level $L$, i.e.
$2^{L_M-L}$, is fixed as a constant.  For a fixed $N$, the time cost
is not monotonic with respect to $L$.  When $L$ is too small the
computational cost of $D^{(L)}$ becomes dominant. When $L$ is too
large, the application of the partial representation $G^{(3)} + \ldots
+ G^{(L)}$ to a vector becomes expensive. From the perspective of time
cost, there is an optimal $L_{\text{opt}}$ for a fixed $N$.  We find
that this optimal level number is the same for $\mc{H}^1$, uniform
$\mc{H}^1$ and $\mc{H}^2$ algorithms.  Table~\ref{tab:time} shows that
$L_{\text{opt}}=4$ for $N=32,64,128$, $L_{\text{opt}}=5$ for $N=256$,
and $L_{\text{opt}}=6$ for $N=512$. This suggests that for large $N$,
the optimal performance is achieved when the size of boxes on the
final level $L$ is $8\times 8$. In other words, $L = L_M - 3$.

The memory cost per DOF for the $\mc{H}^1$, the uniform $\mc{H}^1$ and
the $\mc{H}^2$ algorithms is reported in Table~\ref{tab:memory}.  The
memory cost is estimated by summing the sizes of low-rank
approximations as well as the size of $D^{(L)}$. For a fixed $N$, the
memory cost generally decreases as $L$ increases. This is because as
$L$ increases, an increasing part of the original dense matrix is
represented using low-rank approximations.

\begin{figure}[ht]
  \begin{center}
    \includegraphics[width=2.3in]{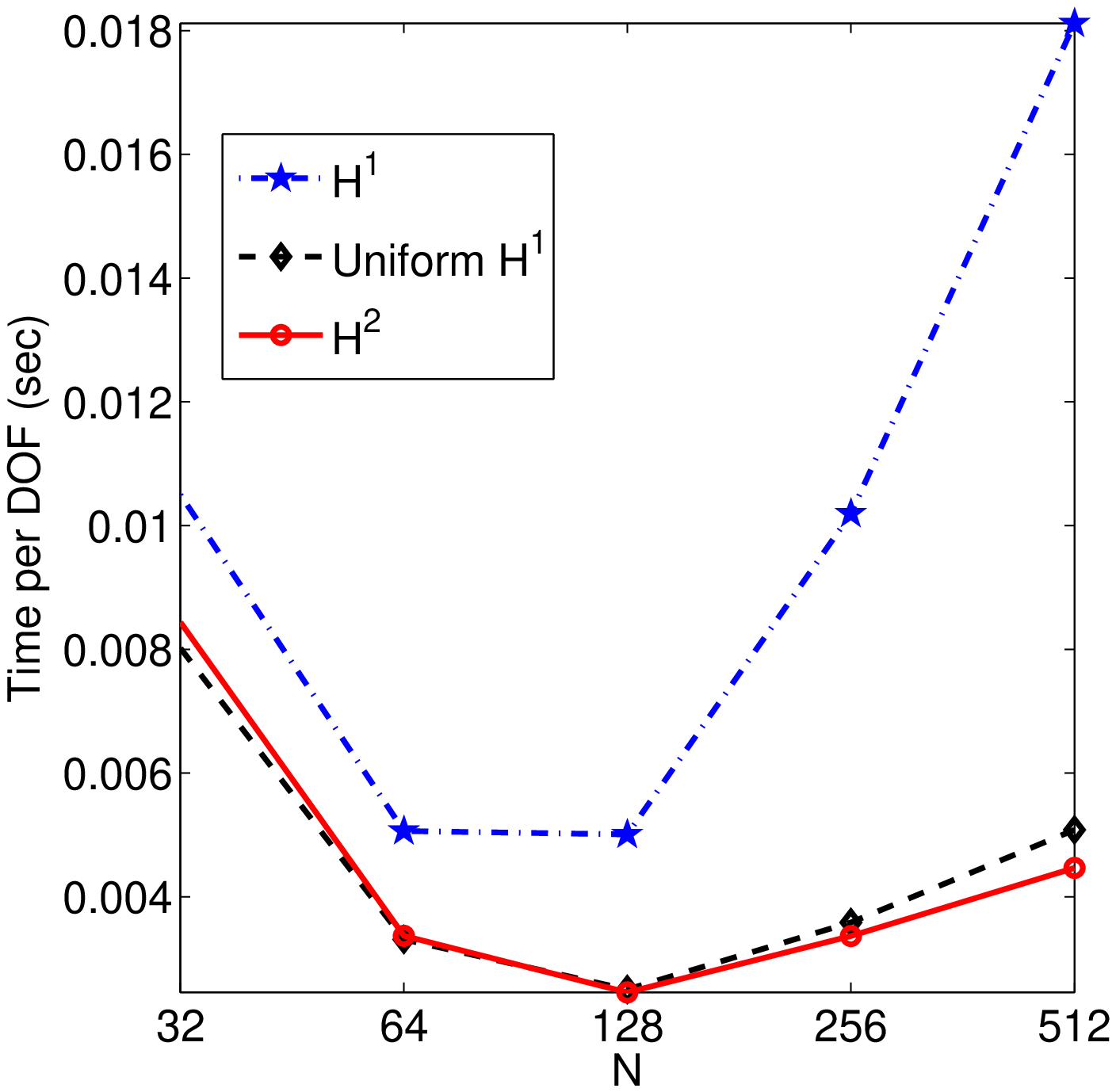}
    \hspace{3em}
    \includegraphics[width=2.2in]{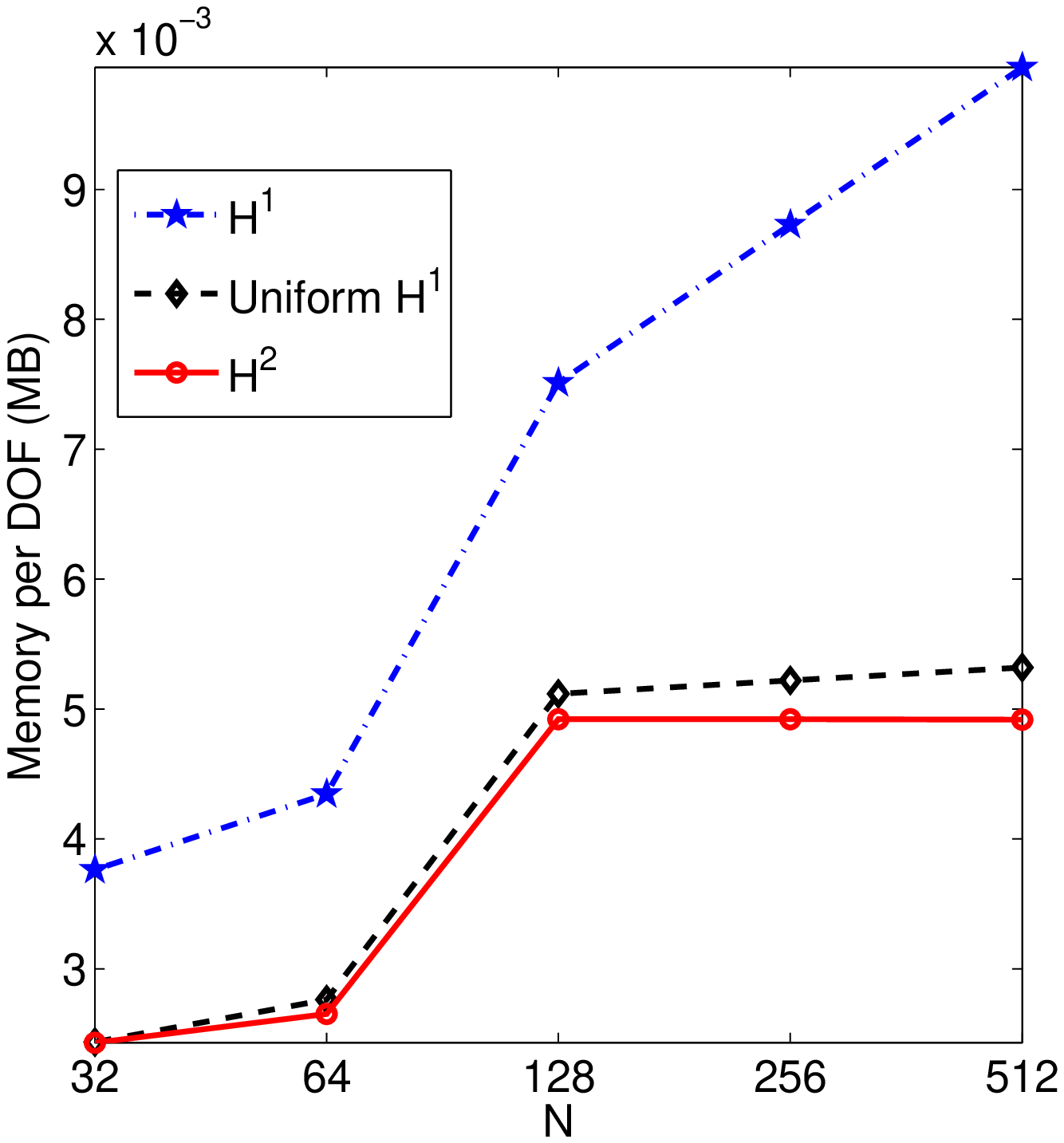}
  \end{center}
  \caption{Comparison of the time and memory costs for the $\mc{H}^1$, the uniform $\mc{H}^1$ and
  the $\mc{H}^2$ versions  with optimal level $L_{\text{opt}}$ for
  $N=32,64,128,256,512$. The x-axis (N) is set to be in logarithmic
  scale.}
  \label{fig:timememory}
\end{figure}

Both Table~\ref{tab:time} and Table~\ref{tab:memory} indicate that
uniform $\mc{H}^1$ algorithm is significantly more advantageous than
$\mc{H}^1$ algorithm, while the $\mc{H}^2$ algorithm leads to a
further improvement over the uniform $\mc{H}^1$ algorithm especially
for large $N$.  This fact can be better seen from
Fig.~\ref{fig:timememory} where the time and memory cost per DOF for
$N=32,64,128,256,512$ with optimal level number $L_{\text{opt}}$ are
shown.  We remark that since the number of calls to the
\textsf{matvec} subroutine are the same in all cases, the time cost
difference comes solely from the efficiency difference of the low rank
matrix-vector multiplication subroutines.

We measure the accuracy for the $\mc{H}^1$, the uniform $\mc{H}^1$ and
the $\mc{H}^2$ representations of $G$ with its actual value using the
operator norm ($2$-norm) of the error matrix. Here, the $2$-norm of a
matrix is numerically estimated by power method~\cite{Golub1996} using
several random initial guesses. We report both absolute and relative
errors. According to Table~\ref{tab:error}, the errors are well
controlled with respect to both increasing $N$ and $L$, in spite of
the more aggressive matrix compression strategy in the uniform
$\mc{H}^1$ and the $\mc{H}^2$ representations.  Moreover, for each box
$\mc{I}$, the rank $r_{U}(\varepsilon)$ of the uniform $\mc{H}^1$
representation is only slightly larger than the rank $r(\varepsilon)$
of the $\mc{H}^1$ representation. This can be seen from
Table~\ref{tab:rankcomparison}. Here the average rank for a level $l$
is estimated by averaging the values of $r_{U}(\varepsilon)$ (or
$r(\varepsilon)$) for all low-rank approximations at level $l$.  Note
that the rank of the $\mc{H}^2$ representation is comparable to or
even lower than the rank in the uniform $\mc{H}^1$ representation.
This is partially due to different weighting choices in the
uniformization step and $\mc{H}^2$ construction step.

\begin{table}[ht]
  \centering
  \begin{tabular}{c|c|c|c|c|c}
    \toprule
    $N$ & $L$ & \textsf{matvec} & $\mc{H}^1$ time   & Uniform $\mc{H}^1$
    time & $\mc{H}^2$ time\\ 
        &   & number          & per DOF (s) & per DOF (s) &per DOF (s) \\
    \midrule
  32 &  4 &      3161 &  0.0106 & 0.0080 &   0.0084  \\ 
    \midrule                                        
  64 &  4 &      3376 &  0.0051 & 0.0033 &   0.0033  \\ 
  64 &  5 &      4471 &  0.0150 & 0.0102 &   0.0106  \\ 
    \midrule                                        
 128 &  4 &      4116 &  0.0050 & 0.0025 &   0.0024  \\ 
 128 &  5 &      4639 &  0.0080 & 0.0045 &   0.0045  \\ 
 128 &  6 &      5730 &  0.0189 & 0.0122 &   0.0125  \\ 
    \midrule                                        
 256 &  4 &      7169 &   0.015 & 0.0054 &   0.0050 \\ 
 256 &  5 &      5407 &   0.010 & 0.0035 &   0.0033 \\ 
 256 &  6 &      5952 &   0.013 & 0.0058 &   0.0057 \\ 
 256 &  7 &      7021 &   0.025 & 0.0152 &   0.0154 \\ 
    \midrule                                        
 512 &  5 &      8439 &   0.025 & 0.0070 &   0.0063   \\ 
 512 &  6 &      6708 &   0.018 & 0.0050 &   0.0044  \\ 
 512 &  7 &      7201 &   0.022 & 0.0079 &   0.0072  \\ 
    \bottomrule
  \end{tabular}
  \caption{\textsf{matvec} numbers and time cost per degree of freedom
  (DOF) for the $\mc{H}^1$, the uniform $\mc{H}^1$ and
  the $\mc{H}^2$ representations with different grid point per dimension $N$ and
  low rank compression level $L$.  The \textsf{matvec} numbers are by
  definition the same in the three algorithms.}
  \label{tab:time}
\end{table}

\begin{table}[ht]
  \centering
  \begin{tabular}{c|c|c|c|c}
    \toprule
    $N$ & $L$ &  $\mc{H}^1$ memory   & Uniform $\mc{H}^1$ memory & $\mc{H}^2$ memory\\ 
        &   &   per DOF (MB) & per DOF (MB) &per DOF (MB) \\
    \midrule
  32 &  4 &        0.0038  & 0.0024 &   0.0024  \\ 
    \midrule                                    
  64 &  4 &        0.0043  & 0.0027 &   0.0026  \\ 
  64 &  5 &        0.0051  & 0.0027 &   0.0026  \\ 
    \midrule                                    
 128 &  4 &        0.0075  & 0.0051 &   0.0049  \\ 
 128 &  5 &        0.0056  & 0.0029 &   0.0027  \\ 
 128 &  6 &        0.0063  & 0.0029 &   0.0027  \\ 
    \midrule                                    
 256 &  4 &         0.0206  & 0.0180 &   0.0177 \\ 
 256 &  5 &         0.0087  & 0.0052 &   0.0049 \\ 
 256 &  6 &         0.0067  & 0.0030 &   0.0027 \\ 
 256 &  7 &         0.0074  & 0.0030 &   0.0027 \\ 
    \midrule                                    
 512 &  5 &         0.0218  & 0.0181 &   0.0177   \\ 
 512 &  6 &         0.0099  & 0.0053 &   0.0049  \\ 
 512 &  7 &         0.0079  & 0.0031 &   0.0027  \\ 
    \bottomrule
  \end{tabular}
  \caption{Memory cost per degree of freedom
  (DOF) for the $\mc{H}^1$, the uniform $\mc{H}^1$ and
  the $\mc{H}^2$ versions with different grid point per dimension $N$ and
  low rank compression level $L$.}
  \label{tab:memory}
\end{table}

\begin{table}[ht]
  \centering
  \begin{tabular}{c|c|c|c|c|c|c|c}
    \toprule
    $N$ & $L$  & \multicolumn{2}{c|}{$\mc{H}^1$} & \multicolumn{2}{c|}{Uniform $\mc{H}^1$} & \multicolumn{2}{c}{$\mc{H}^2$}\\
    \midrule                                                                                                    
    &    & Absolute    & Relative                 & Absolute  & Relative           & Absolute  & Relative          \\
    &    & error & error              & error & error & error & error \\
    \midrule                                                                                                   
    32 &  4      &  2.16e-07 &   3.22e-07         &   2.22e-07 &   3.31e-07       &   2.20e-07 &   3.28e-07       \\ 
    \midrule                                                                                                  
    64 &  4      &  2.10e-07 &   3.15e-07         &   2.31e-07 &   3.47e-07       &   2.31e-07 &   3.46e-07       \\ 
    64 &  5      &  1.96e-07 &   2.95e-07         &   2.07e-07 &   3.12e-07       &   2.07e-07 &   3.11e-07       \\ 
    \midrule                                                                                                  
    128 &  4      &   2.16e-07 &   3.25e-07        &    2.26e-07 &   3.39e-07      &   2.24e-07 &   3.37e-07       \\ 
    128 &  5      &   2.60e-07 &   3.90e-07        &    2.68e-07 &   4.03e-07      &   2.67e-07 &   4.02e-07       \\ 
    128 &  6      &   2.01e-07 &   3.01e-07        &    2.09e-07 &   3.13e-07      &   2.08e-07 &   3.11e-07       \\ 
    \midrule                                                                                                  
    256 &  4      &   1.78e-07 &   2.66e-07        &    1.95e-07 &   2.92e-07      &   2.31e-07 &   3.46e-07       \\ 
    256 &  5      &   2.11e-07 &   3.16e-07        &    2.26e-07 &   3.39e-07      &   2.27e-07 &   3.40e-07       \\ 
    256 &  6      &   2.75e-07 &   4.12e-07        &    2.78e-07 &   4.18e-07      &   2.30e-07 &   3.45e-07       \\ 
    256 &  7      &   1.93e-07 &   2.89e-07        &    2.05e-07 &   3.08e-07      &   2.24e-07 &   3.36e-07       \\ 
    \midrule                                                                                                  
    512 &  5      &   2.23e-07 &   3.35e-07        &    2.33e-07 &   3.50e-07      &   1.42e-07 &   2.13e-07       \\ 
    512 &  6      &   2.06e-07 &   3.09e-07        &    2.17e-07 &   3.26e-07      &   2.03e-07 &   3.05e-07       \\ 
    512 &  7      &   2.67e-07 &   4.01e-07        &    2.74e-07 &   4.11e-07      &   2.43e-07 &   3.65e-07       \\ 
    \bottomrule
  \end{tabular}
\caption{Absolute and relative $2$-norm errors for the $\mc{H}^1$, the uniform $\mc{H}^1$ and
  the $\mc{H}^2$ algorithms with different grid point per dimension $N$ and
  low rank compression level $L$.  The $2$-norm is estimated using power
  method.}
  \label{tab:error}
\end{table}

\begin{table}[ht]
  \centering
  \begin{tabular}{c|c|c|c}
    \toprule
    $l$ & $\mc{H}^1$ & Uniform $\mc{H}^1$ & $\mc{H}^2$ \\
    & average rank   & average rank & average rank\\
    \midrule
    4 & 6 & 13 & 13 \\
    5 & 6 & 13 & 11 \\
    6 & 6 & 12 & 9 \\
    \bottomrule
  \end{tabular}
  \caption{Comparison of the average rank at different levels between
  the $\mc{H}^1$, the uniform $\mc{H}^1$, and the $\mc{H}^2$ algorithms,
  for $N=256$.} 
  \label{tab:rankcomparison}
\end{table}

\FloatBarrier

The peeling algorithm for the construction of hierarchical matrix can
be applied as well to general elliptic operators in divergence form
$H=-\nabla\cdot (a(\bvec{r})\nabla) + V(\bvec{r})$. The computational
domain, the grids are the same as the example above, and five-point
discretization is used for the differential operator. The media is
assumed to be high contrast: $a(i,j) = 1 + U(i,j)$, with $U(i,j)$
being independent random numbers uniformly distributed in $[0,1]$.
The potential functions under consideration are (1)
$V(i,j)=10^{-3}W(i,j)$; (2) $V(i,j)=10^{-6}W(i,j)$.  $W(i,j)$ are
independent random numbers uniformly distributed in $[0,1]$ and are
independent of $U(i,j)$.  We test the $\mc{H}^2$ version for $N=64$,
$L=4$, with the compression criterion $\veps=10^{-6}$. The resulting
$L^2$ absolute and relative error of the Green's function are reported
in Table~\ref{tab:errormedia}. The results indicate that the
algorithms work well in these cases, despite the fact that the
off-diagonal elements of the Green's function have a slower decay than
the first example.  We also remark that the small relative error for
case (2) is due to the large $2$-norm of $H^{-1}$ when $V$ is small.

\begin{table}[ht]
  \centering
  \begin{tabular}{c|c|c}
    \toprule
    Potential &    Absolute error   & Relative error \\
    \midrule                                                                                                   
    $V(i,j)=10^{-3}W(i,j)$ & 5.91e-04 & 2.97e-07\\
    $V(i,j)=10^{-6}W(i,j)$ & 3.60e-03 & 1.81e-09\\
    \bottomrule
  \end{tabular}
  \caption{Absolute and relative $2$-norm errors for the $\mc{H}^2$
  representation of the matrix $\left(-\nabla\cdot (a\nabla) +
  V\right)^{-1}$ with $N=64,L=4$ and two choice of potential function
  $V$.   The $2$-norm is estimated using power method.}
  \label{tab:errormedia}
\end{table}

\section{Conclusions and future work}

In this work, we present a novel algorithm for constructing a
hierarchical matrix from its matrix-vector multiplication. One of the
main motivations is the construction of the inverse matrix of the
stiffness matrix of an elliptic differential operator. The proposed
algorithm utilizes randomized singular value decomposition of low-rank
matrices. The off-diagonal blocks of the hierarchical matrix are
computed through a top-down peeling process. This algorithm is
efficient. For an $n \times n$ matrix, it uses only $\Or(\log n)$
matrix-vector multiplications plus $\Or(n\log n)$ additional steps.
The algorithm is also friendly to parallelization. The resulting hierarchical
matrix representation can be used as a faster algorithm for
matrix-vector multiplications, as well as for numerical homogenization
or upscaling.

The performance of our algorithm is tested using two 2D elliptic operators. The
$\mc{H}^1$, the uniform $\mc{H}^1$ and the $\mc{H}^2$ versions of the
proposed algorithms are
implemented. Numerical results show that our
implementations are efficient and accurate and that the uniform
$\mc{H}^{1}$ representation is significantly more advantageous over $\mc{H}^1$
representation  in terms of both the time cost
and the memory cost, and $\mc{H}^2$ representation leads to further
improvement for large $N$. 

Although the algorithms presented require only $\Or(\log n)$
\textsf{matvec}s, the actual number of \text{matvec}s can be quite
large (for example, several thousands for the example in
Section~\ref{sec:numerical}). Therefore, the algorithms presented here
might not be the right choice for many applications. However, for
computational problems in which one needs to invert the same system
with a huge of unknowns or for homogenization problems where analytic
approaches do not apply, our algorithm does provide an effective
alternative.

The current implementation depends explicitly on the geometric
partition of the rectangular domain. However, the idea of our
algorithm can be applied to general settings. For problems with
unstructured grid, the only modification is to partition the
unstructured grid with a quadtree structure and the algorithms
essentially require no change. For discretizations of the boundary
integral operators, the size of an interaction list is typically much
smaller as many boxes contain no boundary points. Therefore, it is
possible to design a more effective combination strategy with small
number of \textsf{matvec}s. These algorithms can also be extended to
the 3D cases in a straightforward way, however, we expect the constant
to grow significantly. All these cases will be considered in the
future.


\vspace{1em}
\noindent{\bf Acknowledgement:}

L.~L. is partially supported by DOE under Contract
No. DE-FG02-03ER25587 and by ONR under Contract No. N00014-01-1-0674.
L.~Y. is partially supported by an Alfred P. Sloan Research Fellowship
and an NSF CAREER award DMS-0846501. The authors thank Laurent Demanet
for providing computing facility and Ming Gu and Gunnar Martinsson for
helpful discussions. L.~L. and J.~L. also thank Weinan E for support
and encouragement.


\end{document}